\documentclass[12pt,leqno]{amsart} 
\usepackage[numbers]{natbib}
\newcommand{\citeasnoun}{\citet}

\usepackage{amsmath,amsthm,amssymb}
\renewcommand{\widehat}{\hat}
\usepackage{url}
\usepackage[letterpaper,colorlinks={true},
backref,pdfstartview={FitBH},pdfview={FitBH},bookmarks={true}]{hyperref}
\global\hbadness=500
\global\tolerance=200
\def\handlabel{}
\newtheorem*{internallabelledtheorem}{\handlabel}
\newenvironment{labelledtheorem}[1]{%
  \def\handlabel{#1}
\begin{internallabelledtheorem}}{\end{internallabelledtheorem}}

\newtheorem{theorem}{Theorem}[section]
\newtheorem{lemma}[theorem]{Lemma}

\newtheorem{corollary}[theorem]{Corollary} 
\theoremstyle{definition}
\newtheorem{definition}[theorem]{Definition}

\theoremstyle{remark}

\newtheorem{remark}[theorem]{Remark}
\numberwithin{equation}{section}

\newcommand{\ce}{computably enumerable }
\newcommand{\wrt}{w.r.t.\ } 

\newcommand{\join}{\oplus} 
\renewcommand{\phi}{\varphi}

\newcommand{\E}{\mathcal{E}}
\renewcommand{\L}{\mathcal{L}}\renewcommand{\S}{\mathcal{S}}

\newcommand{\R}{\mathcal{R}}
\newcommand{\B}{\mathcal{B}}
\newcommand{\D}{\mathcal{D}}
\newcommand{\C}{\mathcal{C}}

\newcommand{\Ahat}{\widehat{A}} \newcommand{\Bhat}{\widehat{B}}

\newcommand{\Rhat}{\widehat{R}}

\newcommand{\What}{\widehat{W}}\newcommand{\Xhat}{\widehat{X}}

\newcommand{\Shat}{\widehat{S}}\newcommand{\Uhat}{\widehat{U}}
\newcommand{\xhat}{\widehat{x}}\newcommand{\Vhat}{\widehat{V}}

\newcommand{\Cof}{\operatorname{Cof}}
\title[Extensions, Orbits and automorphisms]{Extension theorems,
  Orbits, and automorphisms of the computably enumerable sets}
  
\author[P. Cholak and L. Harrington]{Peter~A.~Cholak and
  Leo~A.~Harrington}

\address{Department of Mathematics\\ University of Notre Dame\\ 
  Notre Dame, IN 46556-5683}

\email{Peter.Cholak.1@nd.edu}

\urladdr{http://www.nd.edu/$\tilde{~}$cholak}

\address{Department of Mathematics\\ University of California \\
  Berkeley, CA 94720-3840}

\email{leo@math.berkeley.edu}

\thanks{Research partially supported NSF Grants DMS-96-34565, 
  99-88716, 02-45167 (Cholak), DMS-96-22290 and DMS-99-71137 (Harrington). We
  would like to thank Bob Soare and Mike Stob for their interest and
  helpful comments.}

\date{August 31, 2005}

\subjclass[2000]{Primary 03D25}

\begin{document}

\begin{abstract} 
  We prove an algebraic extension theorem for the computably
  enumerable sets, $\mathcal{E}$.  Using this extension theorem and
  other work we then show if $A$ and $\widehat{A}$ are automorphic via
  $\Psi$ then they are automorphic via $\Lambda$ where $\Lambda
  \restriction \L^*(A) = \Psi$ and $\Lambda \restriction \E^*(A)$ is
  $\Delta^0_3$.  We give an algebraic description of when an arbitrary
  set $\Ahat$ is in the orbit of a \ce set $A$.  We construct the
  first example of a definable orbit which is not a $\Delta^0_3$
  orbit. We conclude with some results which restrict the ways one can
  increase the complexity of orbits.  For example, we show that if $A$
  is simple and $\widehat{A}$ is in the same orbit as $A$ then they
  are in the same $\Delta^0_6$-orbit and furthermore we provide a
  classification of when two simple sets are in the same orbit.
\end{abstract}

\maketitle

\section{Introduction}

We will work in the structure of the \ce sets.  The language is just
inclusion, $\subseteq$.  This structure is called $\E$. There have
been a large number of papers, see
\cite{Cholak:00,Cholak.Harrington:00,Soare:00} for some recent
surveys, studying $\E$ and the interaction within $\mathcal{E}$ among
the following four mathematical concepts:
\begin{itemize} 
\item Automorphisms: Is there a classification of the orbits of
  $\mathcal{E}$.  Which sets are automorphic, i.e., in the same orbit?
\item Definability: What computably enumerable sets can be defined (in
  the language of just $\{\subset\}$)? Is there a formula which
  distinguishes one set from another within $\mathcal{E}$?
\item Dynamic Properties: How fast (or slow) can a set be enumerated
  compared to another set? or with respect to the standard enumeration
  of all computably enumerable sets?
\item Complexity: How do sets in an orbit interact with each other via
  Turing reducibility? How do the sets in an orbit fit into jump
  classes, in particular, the low$_n$ and high$_n$ classes?  This
  interaction is part of our connection to the \ce degrees.
\end{itemize}
In this paper we focus on automorphisms and orbits although some
aspects of the remaining concepts will arise.

Our understanding of automorphisms of $\E$ is unique to $\E$.  In most
structures with nontrivial automorphisms we can construct
automorphisms via the normal ``back and forth'' argument.  But this is
not the case with $\E$.  To construct automorphisms we use the
properties of being \emph{well-visited} and \emph{well-resided}.
Well-visited is $\Pi^0_2$ and not being well-resided is $\Sigma^0_3$
(we use the negation).  Since the complexity of these properties is at
most $\Sigma^0_3$, the construction of the desired automorphism can be
placed on a tree.  (We will not discuss the details on this placement
nor of the construction of an automorphism of $\E$ but direct the
reader to \citeasnoun{Harrington.Soare:96} or \citeasnoun{Cholak:95}.)
This method is called the $\Delta^0_3$ automorphism method.  If an
automorphism $\Phi$ is constructed on a tree then $\Phi$ has a
presentation computable in the true path (which is $\Delta^0_3$).
Hence all automorphisms constructed in this way are
$\Delta^0_3$-automorphisms. (In some cases we can make the
automorphism effective.)

One step above using the $\Delta^0_3$ automorphism method is to use an
\emph{extension} theorem.  Basically, an extension theorem extends an
isomorphism between two substructures of $\E$ to an automorphism of
$\E$.  The isomorphism between two substructures of $\E$ can be given
in a number of ways and the same can be said about the resulting
automorphism.

Generally, extension theorems are introduced to prove new automorphism
results but they also allow us to reflect back and understand old
automorphism results.  Our philosophy is to argue modularly as much as
possible.  The hope is that an extension theorem provides an
``understandable'' module in the construction of an automorphism of
$\E$.

The first major automorphism result, Soare's result \cite{Soare:74}
that the maximal sets form an orbit, used Soare's Extension Theorem.
In \citeasnoun{Cholak:95}, several more extension theorems were
introduced and used to show that every noncomputable \ce set is
automorphic to a high set.  In \citeasnoun{Cholak:94*1}, the Modified
Extension Theorem was introduced which allowed many of the
automorphism constructions to be recast as using an extension theorem.
For example, in \citeasnoun{Cholak:94*1}, the results about orbits of
hhsimple sets in \citeasnoun{Maass:84} and the result that the
hemimaximal sets form an orbit found in \citeasnoun{Downey.Stob:92}
were recast in this fashion. The Modified Extension Theorem has a
weaker hypothesis than Soare's Extension Theorem.  Soare has recently
proven the ``New Extension Theorem'' and in addition to proving
several new automorphism results with Harrington he has recast almost
all known automorphism results using this and similar theorems (see
\citeasnoun{Soare:00} and \citeasnoun{Soare:nd}).

All of these extension theorems share several common features.  First
they \emph{always} produce $\Delta^0_3$ automorphisms.  All but
Soare's Extension Theorem used the $\Delta^0_3$ automorphism method as
described in \citeasnoun{Cholak:95} and
\citeasnoun{Harrington.Soare:96}. Soare's Extension Theorem was done
effectively.  The isomorphism which these extension theorems extend
and the resulting automorphism are given \emph{dynamically}.

The big issue before applying any extension theorem is to ``match'' up
``entry states'' which is done dynamically. The work done in
Section~\ref{sec:modif-transl-theor} illustrates what we mean by
dynamic, entry states, and matching.

One of the goals of this paper is to prove two new extension theorems
(Theorems~\ref{extension} and \ref{extension2}). These two theorems
differ from the previous extension theorems.  Theorem~\ref{extension2}
allows the possibility that the resulting automorphism is not
$\Delta^0_3$. Both of them are stated ``algebraically'' (or
``statically'').  We have come up with an algebraic description of
entry states and matching using \emph{extendible Boolean algebras} and
\emph{supports}.
Theorem~\ref{extension2} follows algebraically from
Theorem~\ref{extension}.  However we are not free from the use of
dynamic methods.  For example, the proof of Theorem~\ref{extension} is
dynamic and uses Soare's Extension Theorem along with other dynamic
theorems.

(One word of caution: We use the word algebraic to mean facts or
results about the structures we are considering.  The structures we
consider are Boolean algebras and lattices which are ordered
structures where all the definable relations and functions can be
defined just using the order, not necessarily the structures, a model
theorist or algebraist might wish to study.  So a model theorist or
algebraist might wish to read ``order-theoretic'' in place of
``algebraic''.)

Theorem~\ref{extension2} shows that whether an isomorphism between
$\L^*(A)$ and $\L^*(\Ahat)$ can be extended to an automorphism depends
on the existence of a ``nice'' isomorphism among ``some of the entry
states'', where ``some of the entry states'' corresponds to extendible
Boolean algebras and ``nice'' means some properties of the
presentation of the algebras and the isomorphism.

As with any extension theorem, our extension theorems allow us to both
reflect on old automorphism results and prove new automorphism
results.  In Section~\ref{computable}, we reprove some of the
automorphism results mentioned above using Theorems~\ref{extension2}
and \ref{apcomputable}.
One current shortcoming of our extension theorem is with results where
one is given a \ce set $A$ and constructs an automorphic $\Ahat$ with
certain properties (such as highness, for example); this is what Soare
calls a ``type 2'' automorphism result (see
\citeasnoun[Section~7]{Soare:00}). But this might change.

By our extension theorems, the main result from
\citeasnoun{MR2004f:03077} (which depends heavily on
\citeasnoun{MR2003h:03063}) and a result about automorphisms and
extendible Boolean algebras which resembles an automorphism
construction, we can show that if $A$ and $\Ahat$ are automorphic via
$\Psi$ then the isomorphism between $\L^*(A)$ and $\L^*(\Ahat)$
induced via $\Psi$ can be extended into an automorphism $\Lambda$
where $\Lambda \restriction \E^*(A)$ is $\Delta^0_3$.  In other words
we can convert $\Psi$ into an automorphism $\Lambda$ with some nice
properties.

\begin{labelledtheorem}{The Conversion Theorem 
    (Theorem~\ref{mainnew})} If $A$ and $\Ahat$ are automorphic via
  $\Psi$ then they are automorphic via $\Lambda$ where $\Lambda
  \restriction \L^*(A) = \Psi$ and $\Lambda \restriction \E^*(A)$ is
  $\Delta^0_3$.
\end{labelledtheorem}

Hence the complexity of an automorphism comes from the induced
isomorphism between $\L^*(A)$ and $\L^*(\Ahat)$.  The impact of this
theorem is that if we want to show $A$ and $\Ahat$ are automorphic we
are not handicapped by using an extension theorem or the $\Delta^0_3$
automorphism method.  If we show $A$ and $\Ahat$ are automorphic via
$\Lambda$, where $\Lambda$ is built using an extension theorem or the
$\Delta^0_3$ automorphism method, then $\Lambda \restriction \E^*(A)$
is always $\Delta^0_3$.  Our result says if there is an automorphism
taking $A$ to $\Ahat$ then there is an automorphism taking $A$ to
$\Ahat$ which is $\Delta^0_3$ on the inside of $A$ and $\Ahat$.

As a result we get an algebraic description, in terms of the
$\mathcal{L}^*(A)$, $\mathcal{L}^*(\Ahat)$, and extendible algebras,
of when an arbitrary set $\Ahat$ is in the orbit of a \ce set $A$ (see
Theorem~\ref{orbits}).  Not surprisingly the algebraic description is
$\Sigma^1_1$; it begins ``does there exist an isomorphism between
$\L^*(A)$ and $\L^*(\Ahat)$''.

In Section~\ref{sorbits}, we use our extension theorems to show that
there is an elementary definable $\Delta^0_5$ orbit $\mathcal{O}$,
which is not an orbit under $\Delta^0_3$ automorphisms. All the
previously known orbits are orbits under $\Delta^0_3$ automorphisms.

What is surprising is that this complexity comes from how $A \in
\mathcal{O}$ interacts with sets which are disjoint from $A$.  It was
long thought this complexity would come from how $A$ interacts with
sets $W$ such that $W \cap A \neq^* \emptyset$ and $W - A$ is
infinite.  For more details see Section~\ref{restricts2} and
Theorem~\ref{remove}.  In Theorem~\ref{improved}, we improve
Theorem~\ref{remove} to all $A$; we show given an arbitrary \ce set
$A$ the complexity of the orbit of $A$ is determined by the sets
disjoint from $A$.

There will be a sequel to this paper.  In the forthcoming paper we
show that there are orbits which are orbits under $\Delta^0_{\alpha +
  1}$ automorphisms but not $\Delta^0_{\alpha}$ automorphisms, for all
computable $\alpha$.  Cholak, Downey, and Harrington have shown that
the conjecture of Slaman-Woodin that $\{ (A, \Ahat) : A$ is
automorphic to $\Ahat \}$ is $\Sigma^1_1$-complete is correct.  We
hope to use our extension theorems to provide an understandable and
manageable proof of the Slaman-Woodin conjecture.  In fact, we want to
show that there is an $A$ such that whether $\Ahat$ is in the orbit of
$A$ is $\Sigma^1_1$-complete. Theorems~\ref{remove} and \ref{improved}
will have great impact on how we approach these forthcoming results;
they force us to use techniques similar to those used in
Sections~\ref{sec:orbit-mathcalo-4} and
\ref{sec:changes-needed-proof}.  Our extension theorems seem the best
tool for these tasks since we must build non-$\Delta^0_3$
automorphisms in all cases.

Our results certainly justify our philosophy to argue modularly as
much as possible with the use of Soare's Extension Theorem as a
module.  It would be very difficult, if not impossible, to argue that
building automorphisms of $\E$ all at once would be more enlightening.

In Section~\ref{splits}, we introduce and discuss the algebraic
notations needed for our extension theorems.  The remaining sections
have been discussed above.

\section{Splits of $A$}\label{splits}

\subsection{Notation and definitions}

Our notation and definitions are standard and follow
\citeasnoun{Cholak.Harrington:00} which follows \citeasnoun{Soare:87}.

We will be dealing with isomorphisms between various substructures of
$\E$ and automorphisms of $\E$.  In all cases we will think of the
isomorphism (automorphism) as a map from $\omega$ to another copy of
$\omega$, $\widehat{\omega}$.  All subsets of $\widehat{\omega}$ will
wear hats.  We refer to $\widehat{\omega}$ as the \emph{hatted} side
and sometimes we refer to $\omega$ as the \emph{unhatted} side.  When
we define something on the unhatted side there is, of course, the
hatted dual.  We will use this duality frequently without mention.

\subsection{The structure $\S_{\R}(A)$} 

Fix a \ce set $A$.

\begin{definition}\label{badef}
  Let $\mathcal{S}(A)= \{ B : \exists C (B \sqcup C = A)\}$. $\S(A)$
  is the splits of $A$ and $\S(A)$ forms a Boolean algebra.
  $\mathcal{F}(A)$ is the finite subsets of $A$ and is an ideal of
  $\S(A)$.  Let $\mathcal{S}^*(A)$ be the quotient structure
  $\mathcal{S}(A)$ modulo $\mathcal{F}(A)$. Let
  $\mathcal{R}(A) = \{ R : R \subseteq A \text{ and } R %
  \text{ is computable}\}$.  $\R(A)$ is the computable subsets of $A$
  and is an ideal of $\S(A)$.  Let $\mathcal{S}_{\mathcal{R}}(A)$ be
  the quotient structure $\mathcal{S}(A)$ modulo $\mathcal{R}(A)$.
\end{definition}

Let $W$ be in $\S(A)$. Then let $\breve{W} = A - W$ (a \ce set) and
$W^{\R}$ be the equivalence class of $W$ in $\S_{\R}(A)$. From
\citeasnoun[Lemma 2.2]{MR2004f:03077}, we know that if $A$ is
noncomputable, then $\S_\R(A)$ is the atomless Boolean algebra and
hence every Boolean algebra can be embedded in $\S_\R(A)$.

\subsection{$\Sigma^0_3$ Boolean algebras}

Recall from \citeasnoun{Soare:87} the following definition.

\begin{definition}\label{Sigma03BA}
  A countable Boolean algebra $\B = (\{X_i\}_{i \in \omega}, \leq,
  \cup, \cap, \bar{~})$ is a $\Sigma^0_3$ \emph{Boolean algebra} if
  the listing $\{X_i\}_{i \in \omega}$ is uniformly computable and
  there are computable functions $f$ and $g$ and a $\Sigma^0_3$
  relation $R$ such that $X_i \cup X_j = X_{f(i,j)}$, $X_i \cap X_j =
  X_{g(i,j)}$, and $X_i \leq X_j$ iff $R(i,j)$.  (An element of $\B$
  must appear at least once in $\{X_i\}_{i \in \omega}$ but there is
  no bound on the number of times an element may appear in $\{X_i\}_{i
    \in \omega}$.)
\end{definition}

We should be familiar with $\Sigma^0_3$ Boolean algebras.  There is a
beautiful theorem of Lachlan (see \citeasnoun[X.7.2]{Soare:87}) that
says if $\B$ is any $\Sigma^0_3$ Boolean algebra then there is an
hhsimple set $H$ such that $\L^*(H)$ is isomorphic to $\B$.  Let
$\widetilde{\L}(H)$ be the quotient substructure of $\S_{\R}(H)$ given
by $\{R \cap H: R \text{ is computable}\}$ modulo $\R(H)$.  Clearly,
as given, $\widetilde{L}(H)$ is definable in $\E$ with a parameter for
$H$.  In \citeasnoun[Lemma 11.2]{MR2004f:03077}, it is shown that
$\L^*(H)$ and $\tilde{\L}(H)$ are isomorphic.  Hence there is a
substructure of $\S_{\R}(A)$ which ranges over all $\Sigma^0_3$
Boolean algebras as $A$ ranges over all \ce sets.

All of the Boolean algebras we consider will be substructures of
$\S_{\R}(A)$, $\L^*(A)$, or $\E$.  So we will always consider the list
$\{X_i\}_{i \in \omega}$ as a list of \ce sets.  The operations will
be union, intersection, and complementation on \ce sets; and hence the
functions $f$ and $g$ are clearly computable. The relation $R$ will
reflect either $X \subseteq Y$, $X \subseteq_{\R} Y$, or $X
\subseteq^* Y$.

\begin{lemma}
  Given two splits $X$ and $Y$, whether $X \subseteq_{\R} Y$ is
  $\Sigma^0_3$.
\end{lemma}

\begin{proof}
  Given the index for $X$, it is possible to find in a $\Delta^0_3$
  way an index for $\breve{X}$.  Similarly for $Y$. Hence we can find
  an index for $X \triangle Y$ in a $\Delta^0_3$ fashion. Now $X
  \subseteq_{\R} Y$ iff $X \triangle Y$ is computable iff there is an
  $l$ such that $W_l \sqcup (X \triangle Y) = \omega$.  Since ``$W_l
  \sqcup (X \triangle Y) = \omega$'' is $\Pi^0_2$, the last clause in
  the above sentence is $\Sigma^0_3$.
\end{proof}

\begin{theorem}  \label{ba}
  Let $\{X_i : i \in \omega\}$ be a uniformly computable list of \ce
  sets (not necessarily splits of $A$) and a $\Sigma^0_3$ set $B$
  such that $\{X_i : i \in B\}$ generates a subalgebra $\B$ of
  $\S_{\R}(A)$.  Then there is a list $\{Y_i : i \in \omega\}$ where
  all the $Y_i$s are splits of $A$, which witnesses that $\B$ is a
  $\Sigma^0_3$ Boolean algebra.  Furthermore there is a $\Delta^0_3$
  function $g$ from $B$ to $\omega$ such that $X_i = Y_{g(i)}$.
\end{theorem}

\begin{proof}
  Basically we are going to pad the $\Sigma^0_3$ list, $\{X_i : i \in
  B\}$, with lots of finite sets to make it a computable list of \ce
  sets all of which are splits of $A$.  This padding will be done on a
  tree, $2^{< \omega}$.  It will be a standard $\Pi^0_2$ tree
  argument.
  
  Assume $i \in B$ iff $\exists k \phi(i,k)$, where $\phi(i,k)$ is
  $\Pi^0_2$.  Assume that $\phi(i,k)$ is $(\forall x) (\exists
  y)[\Theta(i,k,x,y)]$, where $\Theta$ is $\Delta^0_0$. We define the
  true path by induction as follows: Let $\alpha \subset f$ such that
  $|\alpha| = \langle i , {k} \rangle$.  If $\phi(i,k)$ then $\alpha
  \widehat{~} 0 \subset f$; otherwise $\alpha \widehat{~} 1 \subset
  f$.
  
  The approximation to the true path is also defined by induction.
  Let $\alpha \subseteq f_s$ such that $|\alpha| = \langle i , {k}
  \rangle$ and $|\alpha| \leq s$.  We need a length of agreement
  function: $l_{\alpha}(s)$ is the greatest $z$ such that for all $x
  \leq z$ there is a $y$ with $\Theta(i,k,x,y)$. Let $t < s$ be the
  last stage that $\alpha \subseteq f_s$ (if such a stage does not
  exist let $t=0$).  If $l_{\alpha}(t) < l_{\alpha}(s)$ (an
  $\alpha$-expansionary stage) then $\alpha \widehat{~} 0 \subseteq
  f_s$; otherwise $\alpha \widehat{~} 1 \subseteq f_s$.  It is not too
  hard to show that $f = \liminf_s f_s$.
  
  At $\beta = \alpha \widehat{~} 0$ we will construct a set $Y_j$. If
  $\beta \subseteq f_s$ for the first time ever or the first time
  after being initialized, choose the least $j$ such that $Y_j$ is not
  being constructed and start constructing $Y_j$.  If $\beta \subseteq
  f_s$ and $\beta$ is building $Y_j$, let $Y_{j,s} = X_{i,s}$, where
 $|\alpha| = \langle i , {k} \rangle$.  If
  $\beta$ is to the right of $f_s$ we will initialize $\beta$ at stage
  $s$ (and end the construction of the current $Y_j$).
  
  If $\beta = \alpha \widehat{~} 0 \subset f$ then, by the nature of
  the tree construction, at some stage $\beta$ will be assigned a
  permanent $Y_j$ and never be initialized after that stage. Then $Y_j =
  X_i$, where $|\alpha| = \langle i , {k} \rangle$.  If $Y_j$ is not
  permanently assigned to such a $\beta$ then $Y_j$ is finite.
\end{proof}

\begin{corollary} \label{sraissigma3}
  $\S_\R(A)$ is a $\Sigma^0_3$ Boolean algebra.
\end{corollary}

\begin{proof}
  Given a \ce set $W_e$, it is $\Sigma^0_3$ to decide if $W_e$ is a
  split of $A$ (is there a $j$ such that $W_e \sqcup W_j = A$).
\end{proof}

\begin{definition}
  Following Theorem~\ref{ba}, given $\B$ a $\Sigma^0_3$ Boolean
  algebra of $\S_{\R}(A)$ ($\L^*(A)$ or $\E$), if there is a uniformly
  computable list $\mathcal{X} = \{X_i\}_{i \in \omega}$ of \ce sets
  and a $\Sigma^0_3$ set $B$ such that $\{X_i : i \in B\}$ generates
  $\B$, we say $\mathcal{X}$ and $B$ is a \emph{representation} for
  $\B$. ($B$ might be all of $\omega$.)
\end{definition}

\subsection{Listings of splits of $A$}

We are concerned with the certain well-represented subalgebras of
$\S_{\R}(A)$.  Even if we know $X$ is a split of $A$ we still need
$\mathbf{0''}$ to find a $Y$ such that $X \sqcup Y = A$. We want to
limit ourselves to considering just splits $S$ where we can find $A -
S$ effectively.

\begin{definition}\label{listing}
  A uniformly computable listing, $\S=\{ S_i : i \in \omega \}$, of
  splits of $A$ is an \emph{effective listing} of splits of $A$ iff
  there is another uniformly computable listing $\{ \breve{S}_i : i
  \in \omega \}$ of splits of $A$ such that ${S}_i \sqcup \breve{S}_i
  = A$.
\end{definition}

\begin{lemma}  \label{entrylisting}
  Let $S_e = W_e \searrow A$; this is an entry set. Then the entry
  sets, $\S = \{ S_e: e \in \omega\}$, is an effective listing of
  splits.
\end{lemma}

\begin{proof}
  $(W_e \searrow A) \sqcup (A \backslash W_e) = A$.
\end{proof}

With an entry set the corresponding split is determined at the moment
$x$ enters $A$; either $x$ enters $A$ in $W_e$ or not.  The entry sets
are the canonical example of an effective listing of splits.  This
list depends on the enumeration of $A$.

\begin{lemma}\label{enumeration}
  Let $\S = \{ S_i : i \in \omega \}$ be an effective listing of
  splits of $A$.  Then there is an enumeration of $A$, an effective
  listing of splits of $A$, $\tilde{\S} = \{ \tilde{S}_i : i \in
  \omega \}$, and an effective listing of splits of $A$,
  $\breve{\tilde{\S}} = \{ \breve{\tilde{S}}_i : i \in \omega \}$,
  such that, for all $i$, \wrt the new enumeration of $A$,
  $\tilde{S}_i =^* S_i$, $A \searrow \tilde{S}_i = \emptyset$ (so
  $\tilde{S}_i = \tilde{S}_i \searrow A$), $A \searrow
  \breve{\tilde{S}}_i = \emptyset$,
  $\tilde{S}_i \sqcup \breve{\tilde{S}}_i \sqcup (A \cap \{0,1,\ldots
  i\}) = A$, and if $x \in \tilde{S}_{i,s} \sqcup
  \breve{\tilde{S}}_{i,s}$ then $x \in {S}_{j,s} \sqcup
  \breve{S}_{j,s}$, for all $j \leq i$.
\end{lemma}

\begin{proof}
  Let $x$ enter $A$ (under the old enumeration).  Wait for $x$ to
  enter $S_i$ or $\breve{S}_i$ for $i < x$; adding $x$ to
  $\tilde{S}_i$ or $\breve{\tilde{S}}_i$, respectively.  Then allow
  $x$ to enter $A$ (under the new enumeration).
  
  Clearly $\tilde{\S} = \{ \tilde{S}_i : i \in \omega \}$ and
  $\breve{\tilde{\S}} = \{ \breve{\tilde{S}}_i : i \in \omega \}$ are
  uniformly computable listings of splits of $A$.  The uniformly
  computable listing of splits of $A$, $\{ \breve{\tilde{S}}_i \cup (A
  \cap \{0,1,\ldots, i\}): i \in \omega \}$ witnesses that
  $\tilde{\S}$ is an effective listing of splits.  Similarly $\{
  {\tilde{S}}_i \cup (A \cap \{0,1,\ldots, i\}): i \in \omega \}$
  witnesses that $\breve{\tilde{\S}}$ is an effective listing of
  splits.
\end{proof}

\begin{remark} \label{ness}
  It is necessary that $\mathcal{S}$ be an effective listing of splits
  of $A$ for the above lemma to hold.  The key point of this lemma is
  that when $x$ enters $A$ it has been determined whether $x$ is in
  $\tilde{S}_i$ or not.  So $\tilde{S}_i \sqcup (A \backslash
  \tilde{S}_i) = A$.
  
  \emph{This lemma will be essential.}  It is used in
  Lemma~\ref{entryBA} which in turn plays a key role in
  Section~\ref{sec:meet-hypoth-modif}.  Also see the proof of
  Lemma~\ref{sec:computable}.
\end{remark}

Hence as we vary the enumeration of $A$ we get almost all effective
listing of splits of $A$ as entry sets.  However we do not get all
(noneffective) listing of splits this way.

\begin{lemma}\label{notall}
  No effective listing of splits of infinite \ce set $A$ contains all
  splits of $A$.
\end{lemma}

\begin{proof}
  We will provide two proofs of this lemma.
  
  Let $\S = \{ S_e : e \in \omega \}$ be an effective list of splits
  of $A$.  Let $\{a_i: i \in \omega\}$ be a computable listing of the
  elements of $A$ without repeats.  Let $S = \{ a_i : a_i \notin S_i\}
  = \{ a_i : a_i \in \breve{S}_i \}$.  If $S = S_j$ then $a_j \in S$
  iff $a_j \in S_j$ iff $a_j \notin S_j$. So $S \neq S_j$, for all
  $j$.
  
  By Lemma~\ref{enumeration}, we can assume $S_i = S_i \searrow A$ and
  $\breve{S}_i = \breve{S}_i \searrow A$, for all $i$. By easily
  modifying the Friedberg Splitting Theorem (see
  \citeasnoun[X.2.1]{Soare:87}), we can build a split $S$ and
  $\breve{S}$ such that if $S_i \searrow A$ ($\breve{S}_i \searrow A$)
  is infinite then $S_i \searrow S$ ($\breve{S}_i \searrow S$) is
  infinite and similarly for $\breve{S}$.  The split $S$ is not in
  $\S$.
\end{proof}

\subsection{Extendible subalgebras}

We would like to consider subalgebras of $\S_{\R}(A)$ which have a
representation that is an effective listing of splits of $A$.   

\begin{definition}
  A $\Sigma^0_3$ subalgebra $\B$ of $\S_{\R}(A)$ is \emph{extendible}
  iff there is representation $\mathcal{S}$ and $B$ of $\B$ such that
  $\mathcal{S}$ is an effective listing of splits of $A$ and $B$
  is a $\Delta^0_3$ set.
\end{definition}

We will assume that if $\B$ is extendible then the given
representation is always an effective listing of splits of $A$. From
this point further $\S =\{ S_i: i \in \omega\}$ will always refer to
an effective listing of splits of $A$ and $\mathcal{X} = \{ X_i: i \in
\omega\}$ to a uniformly computable list of \ce sets.

\begin{lemma}\label{trivial}
  The trivial subalgebra of $\S_\R(A)$ is extendible.
\end{lemma}

\begin{proof}
  Let $S_{2e} = \emptyset$, $\breve{S}_{2e} = A$, $S_{2e+1} = A$,
  $\breve{S}_{2e+1} = \emptyset$, and $B = \omega$.
\end{proof}

\begin{lemma}\label{entry}
  The subalgebra $\E_A$ generated by the entry sets is extendible
  (this is what we call an entry set Boolean algebra for $A$).
\end{lemma}

\begin{proof}
   Use the listing from Lemma~\ref{entrylisting} and $B = \omega$. 
\end{proof}

\begin{lemma}\label{entryBA}%
  Let $\B \subseteq \S_{\R}(A)$ be extendible via $\mathcal{S}$ and
  $B$.  There is an enumeration of $A$ and an effective listing of
  splits, $\tilde{\S} = \{ \tilde{S}_i : i \in \omega\}$, such that
  $\tilde{\S}$ and $B$ witness that ${\B}$ is extendible and, for
  all $i$, $A \searrow \tilde{S}_i = \emptyset$ (and so $\tilde{S}_i
  \sqcup (A \backslash \tilde{S}_i) = A$).
\end{lemma}

\begin{proof}%
  Apply Lemma~\ref{enumeration} to $\mathcal{S}$ to get the desired
  enumeration of $A$ and the effective listing of splits of $A$,
  $\tilde{\S}$. $\{\tilde{S}_i : i \in B\}$ generates $\B$.
\end{proof}

Hence every extendible Boolean algebra is an extendible subalgebra of
an entry set Boolean algebra.  Clearly every extendible Boolean 
algebra is a $\Sigma^0_3$ Boolean
algebra.

\begin{lemma}\label{joinext}
  If $\B$ and $\B'$ are extendible then $\B \join \B'$ are
  extendible.  
\end{lemma}

\begin{proof}
  Let $\{S_i\}_{i \in \omega}$ and $B$ witness that $\B$ is extendible
  and similarly for $\B'$.  Let $T_{2i} = S_i$ and $T_{2i+1}= S'_i$.
  Then $\{T\}_{i \leq \omega}$ and $\{2i : i \in B\} \cup \{2i+1: i
  \in B'\}$ witness that $\B \join \B'$ is extendible.
\end{proof}

\begin{theorem}\label{sec:extend-subalg}
  There is an extendible algebra $\B$ of $\S_\R(A)$ such that
  \begin{enumerate}
  \item for all $i \in B$, $S_i$ is computable,
  \item for all $R \in \R(A)$, there is $i \in B$ such that $R = S_i$,
    and
  \item $B$ is infinite.
  \end{enumerate}
\end{theorem}

\begin{proof}
  For this proof fix an enumeration of $A$ (with $A_1 = \emptyset$).
  The idea is that if $R$ is a computable split of $A$ then there are
  $i_0, i_1, i_2 $ such that $R = W_{i_0}$, $A \searrow W_{i_0} =
  \emptyset$ (\wrt this fixed enumeration), $\breve{W}_{i_0} =
  W_{i_1}$, $A \searrow W_{i_1} = \emptyset$, $W_{i_0,s+1} \sqcup
  W_{i_1,s+1} = A_{s+2}$, $\overline{W}_{i_0}= W_{i_2}$, and
  $W_{i_1,s+1} \subseteq W_{i_2,s+1}$, for all $s$, (before $x$ enters
  $A$ determine which of $R$ or $\overline{R}= W_{i_2}$ $x$ is in and
  add $x$ to $W_{i_0}$ or $W_{i_1}$ and $W_{i_2}$ accordingly). In
  this case, we can let $S_i = W_{i_0}$ and $\breve{S}_i = W_{i_1}$,
  where $i = \langle i_0, i_1, i_2 \rangle$.  But to make $\S$ a
  uniformly computable list of \ce sets we must be more careful.

  Let $i = \langle i_0, i_1, i_2 \rangle$.  Assume that $S_{i,s}$ and
  $\breve{S}_{i,s}$ have been defined and $i$ has not been declared
  \emph{unusable}.  If $(A \searrow W_{i_0})_{s+1} = \emptyset$, $(A
  \searrow W_{i_1})_{s+1} = \emptyset$, $W_{i_0,s+1} \sqcup
  W_{i_1,s+1} = A_{s+2}$, $W_{i_0,s+1} \cap W_{i_2,s+1} = \emptyset$,
  and $W_{i_1,s+1} \subseteq W_{i_2,s+1}$, then let $S_{i,s+1} =
  W_{i_0,s+1}$ and $\breve{S}_{i,s+1} = W_{i_1,s+1}$.  Otherwise
  declare $i$ \emph{unusable} and, for all $s' > s$, let $S_{i,s'} =
  S_{i,s}$ and $\breve{S}_{i,s'} = A_{s'+1} - S_{i,s}$.  $\{S_i\}_{i
    \in \omega}$ is an effective listing of splits of $A$.
  
  Let $i \in B$ iff $W_{i_0} \sqcup W_{i_1} = A$, $A \searrow W_{i_0}
  = \emptyset$, $A \searrow W_{i_1} = \emptyset$, $W_{i_0,s+1} \sqcup
  W_{i_1,s+1} = A_{s+2}$, $\overline{W}_{i_0}= W_{i_2}$, and
  $W_{i_1,s+1} \subseteq W_{i_2,s+1}$, for all $s$.  $B$ is
  $\Delta^0_3$.
  
  $\{S_i\}_{i \in \omega}$ and $B$ represent our extendible algebra
  $\B$.  If $i \in B$ then $S_i = W_{i_0}$, $\breve{S}_i = W_{i_1}$,
  and $S_i \sqcup W_{i_2} = \omega$ and hence $S_i$ is computable.
  Given a computable subset $R$ of $A$, by the first paragraph of this
  proof, there is an corresponding $i \in B$ with $R = W_{i_0}$.
  Since there are infinitely many such $R$, $B$ is infinite.
\end{proof}

\subsection{Isomorphisms}

\begin{definition}\label{defiso}
  We consider $\Theta$ a partial map between splits of $A$ and splits
  of $\Ahat$ an \emph{isomorphism} between a substructure $\B$ of
  $\S_\R(A)$ and a substructure $\widehat{\B}$ of $\S_\R(\Ahat)$ if
  $\Theta$ preserves $\subseteq_\R$, for each equivalence class
  $S_{\R}$ of $\B$ if $S \in S_\R$, $\Theta(S)$ exists, and for each
  equivalence class $\widehat{S}_{\R}$ of $\widehat{\B}$ if $\Shat \in
  \widehat{S}_\R$, $\Theta^{-1}(\Shat)$ exists. There is a function
  $h$ such that $\Theta(W_i) = \widehat{W}_{h(i)}$ and
  $\Theta^{-1}(\widehat{W}_i) = W_{h^{-1}(i)}$.  If $h$ is
  $\Delta^0_3$ then so is $\Theta$. 
\end{definition}

\begin{definition}\label{extiso}
  We say two extendible Boolean algebras $\B$ and $\widehat{\B}$ are
  \emph{extendibly isomorphic} via $\Theta$ iff 
  \begin{itemize}
  \item there is an effective listing of splits $\{S_i\}_{i \in
      \omega}$ and a $B$ which witness that $\B$ is an extendible
    algebra,
  \item there are $\{\Shat_i\}_{i \in \omega}$ and $\widehat{B}$ which
    witness $\widehat{\B}$ is an extendible algebra,
  \item for all $i \in B$, there is a $j \in \widehat{B}$ such that
    $\Theta(S_i) = \Shat_j$,
  \item for all $j \in \widehat{B}$ there is an $i \in {B}$ such that
    $\Theta^{-1}(\Shat_j) = S_i$, and
  \item this partial map induces an isomorphism $\Theta'$ between $\B$
    and $\widehat{\B}$ as in Definition~\ref{defiso}.
 \end{itemize}
   In this case, we say that $\Theta$ is an
  \emph{extendible isomorphism}.  There is a function $h$ such that
  $\Theta(S_i) = \widehat{S}_{h(i)}$ and $\Theta^{-1}(\widehat{S}_i) =
  S_{h^{-1}(i)}$.  If $h$ is $\Delta^0_3$ then so is $\Theta$.  We
  write $\Theta(S_i) = \Shat_{\Theta(i)}$ and
  $\Theta^{-1}(\Shat_j)=S_{\Theta^{-1}(j)}$. If $S$ is not an $S_i$, for
  all $i$, but $S_\R \in \B$ we let $\Theta(S) = \Theta'(S)$ and
  similarly for $\Shat$. Hence we will also consider $\Theta$ to be an
  isomorphism (as in Definition~\ref{defiso}) between $\B$ and
  $\widehat{\B}$.
\end{definition}

\begin{lemma}\label{sec:preserve}
  Let $\B$ be a $\Sigma^0_3$ substructure of $\S_\R(A)$ and
  $\widehat{\B}$ be a $\Sigma^0_3$ substructure of $\S_\R(\Ahat)$.
  Assume that $\Theta$ is a map between $\{X_i : i \in B\}$ and
  $\{\widehat{X}_i : i \in \widehat{B}\}$.  Furthermore assume that
  for $i,j \in B$, $X_i-X_j$ is computable iff
  $\Theta(X_i)-\Theta(X_j)$ is computable and, dually, for all $i,j
  \in \widehat{B}$, $\Xhat_i-\Xhat_j$ is computable iff
  $\Theta^{-1}(\Xhat_i)-\Theta^{-1}(\Xhat_j)$ is computable.  Then
  $\Theta$ induces an isomorphism $\Theta'$ between $\B$ and $\widehat{\B}$.
\end{lemma}

\begin{proof}
  $\Theta$ and $\Theta^{-1}$ preserve $\subseteq_\R$. $X_j
  \subseteq_\R X_i$ iff $X_j - X_i$ is computable iff $\Theta(X_j) -
  \Theta(X_i)$ is computable iff $\Theta(X_j) \subseteq_\R
  \Theta(X_i)$. And similarly for $\Theta^{-1}$.  Given $S_\R \in \B$
  find $i$ such that $X_i \in S_\R$ and, for all $S \in S_\R$, let
  $\Theta'(S) = \Theta(X_i)$.  $\Theta'$ is well defined and preserves
  $\subseteq_\R$ since $\Theta$ does. Define $\Theta^{-1}$ dually.
\end{proof}

If $\Theta$ is an extendible isomorphism and we apply
Lemma~\ref{entryBA} to the effective listing of splits then $\Theta$
remains an extendible isomorphism between these two extendible
algebras with regard to the new listing of splits. 

\begin{lemma}\label{trivialiso}
  The trivial subalgebras of $\S_\R(A)$ and $\S_\R(\Ahat)$ are
  effectively extendibly isomorphic as extendible subalgebras of
  $\S_\R(A)$ and $\S_\R(\Ahat)$.
\end{lemma}

\begin{proof}
  Let $\{S_i\}_{i < \omega}$ be the listing of splits given in
  Lemma~\ref{trivial} for the trivial subalgebra of $\S_\R(A)$.
  Let $\{\Shat_i\}_{i < \omega}$ be the listing of splits given in
  Lemma~\ref{trivial} for the trivial subalgebra of $\S_\R(\Ahat)$.
  Let $\Theta(S_i) = \Shat_i$ and $\Theta^{-1}(\Shat_i) = S_i$.
\end{proof}

\begin{lemma}\label{joinextiso}
  Assume that $\B$ and $\widehat{\B}$ are extendible subalgebras which are
  extendibly isomorphic via $\Theta$. Assume that ${\B}'$ and
  $\widehat{\B}'$ are extendible subalgebras which are extendibly
  isomorphic via $\Theta'$.  Then, by Lemma~\ref{joinext}, $\B \join
  \B'$ and $\widehat{\B} \join \widehat{\B}'$ are extendible
  subalgebras which are extendibly isomorphic via $\Delta$, where
  $\Delta(T_{2e}) = \Theta(S_e)$, $\Delta(T_{2e+1}) = \Theta'(S'_e)$,
  $\Delta^{-1}(\widehat{T}_{2e}) = \Theta^{-1}(\widehat{S}_e)$, and
  $\Delta^{-1}(\widehat{T}_{2e+1}) = (\Theta')^{-1}(\widehat{S}'_e)$.
\end{lemma}

\section{Extensions to isomorphisms}\label{ext}

Recall that $\E^*(A)$ is the structure $(\{W_e \cap A : e \in
\omega\}, \subseteq)$ modulo the finite sets.  An isomorphism between
$\E^*(A)$ and $\E^*(\Ahat)$ is a one-to-one, onto (both of these items
are in terms of $*$-equivalence classes) function, $\Xi$, from $\{W_e
\cap A : e \in \omega\}$ to $\{\What_e \cap \Ahat: e \in \omega\}$
such that $W_e \cap A \subseteq^* W_i \cap A$ iff $\Xi(W_e \cap \Ahat)
\subseteq^* \Xi(W_i \cap \Ahat)$.  Note the $\Xi$ is applied to $W_e
\cap A$, not $W_e$.

The goal of this section is to prove and discuss the import of the
following extension theorem.  

\begin{theorem}\label{extension}
  Let $\B \subseteq \S_{\R}(A)$ and $\widehat{\B} \subseteq
  \S_{\R}(\Ahat)$ be two extendible Boolean algebras which are
  $\Delta^0_3$ extendibly isomorphic via $\Theta$.  Then there is a
  $\Phi$ such that $\Phi$ is a $\Delta^0_3$ isomorphism between
  $\E^*(A)$ and $\E^*(\Ahat)$, for all $i \in B$, $\Phi(S_i) =_{\R}
  \Theta(S_i)$, and for all $i \in \widehat{B}$, $\Phi^{-1}(\Shat_i)
  =_{\R} \Theta^{-1}(\Shat_i)$.
\end{theorem}

What is important about this theorem is that we can \emph{extend} the
extendible isomorphism between $\B$ and $\widehat{\B}$ to an
isomorphism between $\E^*(A)$ and $\E^*(\Ahat)$.

The first clause of the conclusion should not be very surprising.
After all, if $A$ and $\Ahat$ are infinite then there is an effective
isomorphism $\Psi$ between $\E^*(A)$ and $\E^*(\Ahat)$. Let $f$ be an
effective map from $A$ to $\Ahat$ and $\Psi(W) = f(W)$.  Moreover, if
$A$ and $\Ahat$ are computable then $\Psi$ clearly computably agrees
with $\Theta$ on all $S_i$ and hence the second clause of the
conclusion holds with $\Psi$. 

The main use of Theorem~\ref{extension} is in the proof of
Theorem~\ref{extension2} and Theorem~\ref{extension4}.  These are the
only examples of the use of Theorem~\ref{extension} in this paper.
However, we will provide several examples of the use of
Theorem~\ref{extension2} and Theorem~\ref{extension4}.

There are several possible ways to prove this theorem. For example,
one could use some of Soare's recent work on extension theorems.  We
had used such a proof in an earlier version of this paper.  In this
version we will base our proof on published theorems. However, we will
have to use them in novel ways and, in a few cases, note that these
proofs prove more than what is actually stated.

We will base our proof on a theorem, the Translation Theorem, from
\citet{Cholak:94*1}.  The proof will have a few parts.  First we will
restate the Translation Theorem in a slightly strengthened form and
show why this version follows from the proof in \citet{Cholak:94*1}.
Then we construct a $\bf{0''}$ enumeration witnessing that $\Theta$
is an extendible isomorphism and meeting the hypothesis of the
Translation Theorem.  Then we apply the modified Translation Theorem
followed by Soare's original Extension Theorem to this enumeration to
get the desired isomorphism.

The proof of Theorem~\ref{extension} is one of the few places where we
have to go into the difficult details of actually building an
isomorphism by a dynamic construction and the use of states.

\subsection{The Modified Translation Theorem}
\label{sec:modif-transl-theor}

These next definitions are a repeat of the first six definitions in
Section~1 of \citet{Cholak:94*1} using slightly different notation.

\begin{definition}
   \begin{enumerate}
   \item $\{X_n\}_{n < \omega}$ is a \emph{uniformly computable
       collection} of c.e.\ sets if there is a computable function $h$
     such that for all $n$, $X_n = W_{h(n)}$.
     
   \item $\{X_n\}_{n < \omega}$ is a \emph{uniformly
       $\bf{0''}$-computable collection} of c.e.\ sets if there is a
     function $h \leq_T \bf{0''}$ such that for all $n$, $X_n =
     W_{h(n)}$.
     
   \item $\{X_{n,s}\}_{n < \omega, s < \omega}$ is a \emph{uniformly
       $\bf{0''}$-computable enumeration} of c.e.\ sets if there is a
     function $h \leq_T \bf{0''}$ such that for all $n$ and $s$,
     $X_{n,s} = W_{h(n),s}$.
   \end{enumerate}
\end{definition}

\begin{definition}
  For any $e$, if we are given uniformly computable enumerations of
  $\{X_{n,s}\}_{n\leq e, s < \omega}$ and $\{Y_{n,s}\}_{n\leq e, s <
    \omega}$ of c.e.\ sets $\{X_n\}_{n \leq e}$ and $\{Y_n\}_{n \leq
    e}$, define the \emph{full $e$-state of $x$ at stage $s$,
    $\nu(e,x,s)$}, with respect to (w.r.t.)  $\{X_{n,s}\}_{n \leq e,s
    < \omega}$ and $\{Y_{n,s}\}_{n \leq e,s < \omega}$ to be the
  triple
  \begin{equation*}
    \nu(e,x,s) = \langle e ,\sigma(e,x,s), \tau(e,x,s) \rangle
  \end{equation*}
   where
   \begin{equation*}
      \sigma(e,x,s) = \{i \leq e: x \in X_{i,s}\}
   \end{equation*}
    and

   \begin{equation*}
      \tau(e,x,s) = \{i \leq e: x \in Y_{i,s}\}.
   \end{equation*}
\end{definition}

\begin{definition}
  For any collection of c.e.\ sets $\{X_n\}_{n \leq e}$ and
  $\{Y_n\}_{n \leq e}$, define the \emph{final $e$-state of $x$,
    $\nu(e,x)$}, w.r.t\ $\{X_n\}_{n \leq e}$ and
  $\{Y_n\}_{n \leq e}$ to be the triple
  \begin{equation*}
    \nu(e,x) = \langle e ,\sigma(e,x), \tau(e,x) \rangle
  \end{equation*}
   where
   \begin{equation*}
      \sigma(e,x) = \{i \leq e: x \in X_{i}\}
   \end{equation*}
    and
   \begin{equation*} 
      \tau(e,x) = \{i \leq e: x \in Y_{i}\}.
   \end{equation*}
\end{definition}

\begin{definition}\label{entrystate}
  Assume that $\{A_s\}_{s <\omega}$ is a uniformly computable
  enumeration of $A$, an infinite c.e.\ set.  For any $e$, assume we
  are given uniformly computable enumerations of $\{X_{n,s}\}_{n\leq
    e, s < \omega}$ and $\{Y_{n,s}\}_{n\leq e, s < \omega}$ of c.e.\ 
  sets $\{X_n\}_{n \leq e}$ and $\{Y_n\}_{n \leq e}$.  For each full
  $e$-state $\nu$, define the c.e.\ set
  \begin{multline*}
    D^A_\nu = \{x : \exists t \text{ such that } x \in A_{s+1} - A_{s}
    \text{ and } \nu = \nu(e,x,s) \\ \text{ w.r.t. }
    \{X_{n,s}\}_{n\leq e, s < \omega} \text{ and } \{Y_{n,s}\}_{n\leq
      e, s < \omega}\}.
  \end{multline*}
  If $x \in D^A_\nu$, we say that $\nu$ is the \emph{entry state} of
  $x$ w.r.t.\ $\{X_{n,s}\}_{n\leq e, s < \omega}$ and
  $\{Y_{n,s}\}_{n\leq e, s < \omega}$ into $A$. We say that $D^A_\nu$
  is measured w.r.t.\ $\{X_{n,s}\}_{n\leq e, s < \omega}$ and
  $\{Y_{n,s}\}_{n\leq e, s < \omega}$.
\end{definition}

The following definition is new and is used for notation ease.

\begin{definition}
  We write $X \doteq_\R Y$ iff $X \subseteq Y$ and $X =_\R Y$.
\end{definition}

\begin{theorem}[\bf{The Modified Translation Theorem}] 
  \label{mtt}
  Assume that $\{A^\dagger_s\}_{s \in \omega}$,
  $\{\Ahat^\dagger_s\}_{s \in \omega}$, $\{U^\dagger_{n,s}\}_{n <
    \omega, s < \omega}$, $\{\Vhat^\dagger_{n,s}\}_{n < \omega, s <
    \omega}$, $\{\Uhat^\dagger_{n,s}\}_{n< \omega, s < \omega}$, and
  $\{V^\dagger_{n,s}\}_{n < \omega, s < \omega}$ are uniformly
  $\bf{0''}$-computable enumerations of the infinite c.e.\ sets
  $A^\dagger $ and $\Ahat^\dagger $ and the uniformly
  $\mathbf{0''}$-computable collection of c.e.\ sets
  $\{U^\dagger_{n}\}_{n< \omega}$, $\{\Vhat^\dagger_{n}\}_{n <
    \omega}$, $\{\Uhat^\dagger_{n}\}_{n< \omega}$, and
  $\{V^\dagger_{n}\}_{n< \omega}$ satisfying the following conditions:
  \begin{equation} \label{eq:1}
        (\forall n) [ \Ahat^\dagger \searrow \Uhat^\dagger_n
        =  A^\dagger
        \searrow \Vhat^\dagger_n = \emptyset], 
  \end{equation}
  \begin{equation} \label{eq:2}
        (\forall \nu)[D^{A^\dagger}_\nu \text{ is infinite iff
        }D^{\Ahat^\dagger}_{\nu} \text{ is
        infinite}],
  \end{equation}
  where, for all $e$-states, $D^{A^\dagger}_\nu$ is measured w.r.t\ 
  $\{U^\dagger_{n,s}\}_{n\leq e, s < \omega}$ and
  $\{\Vhat^\dagger_{n,s}\}_{n\leq e, s < \omega}$, and
  $D^{\Ahat^\dagger}_\nu$ is measured w.r.t\ 
  $\{\Uhat^\dagger_{n,s}\}_{n\leq e, s < \omega}$ and
  $\{V^\dagger_{n,s}\}_{n\leq e, s < \omega}$.
  
  Then there is a collection of uniformly computable c.e.\ sets
  $\{U_{n}\}_{n< \omega}$, $\{\Vhat^{+}_{n}\}_{n< \omega}$,
  $\{\Uhat^{+}_{n}\}_{n< \omega}$, and $\{V_{n}\}_{n< \omega}$ and
  uniformly computable enumerations $\{A_s\}_{s \in \omega}$,
  $\{\Ahat_s\}_{s \in \omega}$, $\{U_{n,s}\}_{n< \omega, s < \omega}$,
  $\{\Vhat^{+}_{n,s}\}_{n< \omega, s < \omega}$,
  $\{\Uhat^{+}_{n,s}\}_{n< \omega, s < \omega}$, and $\{V_{n,s}\}_{n <
    \omega, s < \omega}$ of these sets such that
  \begin{equation}
    \label{eq:3}
    A_{s+1} = A^\dagger_s \text{ and } \Ahat_{s+1} = \Ahat^\dagger_s,
  \end{equation}
  \begin{equation}
    \label{eq:4}
    (\forall n) [ \Ahat \searrow \Uhat^{+}_n = 
    A \searrow \Vhat^{+}_n =\emptyset],
  \end{equation}
  \begin{equation}
    \label{eq:5}
    (\forall n) (\exists e_n) [ U^\dagger_n =^* U_{e_n},
     \Vhat^{+}_{e_n}  \doteq_\R \Vhat^\dagger_n,
    \Uhat^{+}_{e_n} \doteq_\R \Uhat^\dagger_n, 
     \text{ and } V^\dagger_n =^* V_{e_n}], 
   \end{equation}
   \begin{equation}
     \label{eq:6}
     \begin{split}
       (\forall e)[\text{either }[U_e \backslash A =^* \Vhat^{+}_e
       \backslash A =^* \Uhat^{+}_e \backslash \Ahat =^* V\backslash
       \Ahat =^* \emptyset] \\(\text{hence, by Equation~\eqref{eq:4}, }
       \Uhat^{+}_e = \Vhat^{+}_e =^* \emptyset) \text{ or } \\ [
       \text{there is an } n \text{ such that } e = e_n ~ (\text{from
         Equation~\eqref{eq:5}}) ]],
     \end{split}
   \end{equation}
   \begin{equation} \label{eq:7}
        (\forall \nu)[D^{\Ahat}_\nu \text{ is infinite implies
        }(\exists \nu' \geq \nu)D^A_{\nu'} \text{ is
        infinite}],
  \end{equation}
   \begin{equation} \label{eq:8}
        (\forall \nu)[D^{A}_\nu \text{ is infinite implies
        }(\exists \nu' \leq \nu)[D^{\Ahat}_{\nu'} \text{ is
        infinite}]],
  \end{equation}
  where, for all $e$-states, $D^{A}_\nu$ is measured w.r.t\ 
  $\{U_{n,s}\}_{n\leq e, s < \omega}$ and $\{\Vhat^{+}_{n,s}\}_{n\leq
    e, s < \omega}$, and $D^{\Ahat}_\nu$ is measured w.r.t\ 
  $\{\Uhat^{+}_{n,s}\}_{n\leq e, s < \omega}$ and $\{V_{n,s}\}_{n\leq
    e, s < \omega}$.
\end{theorem}

\subsection{Proving the Modified Translation Theorem} 

We will show that the Modified Translation Theorem follows from the version
of the Translation Theorem published in \citet{Cholak:94*1}.
Equations labeled ``3.x'' refer to the Modified Translation Theorem
and equations labeled ``1.x'' refer to the Translation Theorem.

First note that rather than $A^\dagger$, $A$, $\Ahat^\dagger$,
$\Ahat$, $\Uhat^{+}$, and $\Vhat^{+}$ the published version of the
Translation Theorem used $T^\dagger$, $T$, $\widehat{T}^\dagger$,
$\widehat{T}$, $\Uhat$, and $\Vhat$.  So Equation~\ref{eq:1} is the
same as Equation~1.7.  Equation~\ref{eq:2} implies Equations~1.8 and
1.9.  Hence this version is weaker than the published version.  We
could weaken the hypothesis of this version but for our current uses
there is no need.

In the conclusions, Equation~\ref{eq:3} is the same as Equation~1.10,
Equation~\ref{eq:4} is the same as Equation~1.11, Equation~\ref{eq:7}
is the same as Equation~1.14, and Equation~\ref{eq:7} is the same as
Equation~1.15.

That leaves Equations~\ref{eq:5} and \ref{eq:6}.  Equations~1.12 and
1.13 are shown true on page~95 of \citet{Cholak:94*1} (lines -13 to
-11). (Note in Equation~1.12, the first and only ``$\cup$'' should be
a ``$\cap$''.)  We will start from the middle of page~95 and show that
Equations~\eqref{eq:5} and \eqref{eq:6} hold.

Recall $g$ is an onto, one-to-one, computable function from $\omega$ to
$Tr$. In \cite{Cholak:94*1}, $U_e = U_{g(e)}$ and similarly for
$\Vhat^{+}$, $\Uhat^{+}$, and $V$, while $U^\dagger_{g(e)} =
U^\dagger_{|g(e)|}$ and similarly for $\Vhat^\dagger$, $\Uhat^\dagger$,
and $V^\dagger$. If $g(e) \not\subset f$ then the first clause of
Equation~\eqref{eq:6} holds.  If $\beta = g(e) \subset f$ and $n =
|g(e)|$ then it is enough to show $e=e_n$.  (That is, it is enough to
show Equation~\eqref{eq:5} holds for $n$ and $e$.)  So rather than
showing $\Vhat^\dagger_n \cap \overline{A} =^* \Vhat^{+}_{e} \cap
\overline{A}$ we must show $\Vhat^{+}_{e} \doteq_\R \Vhat^\dagger_n $
and similarly for $\Uhat^{+}$ and $\Ahat$ and we will be done.

By Lemma~2.12 of \citet{Cholak:94*1}, the fact that for all $x$,
$\alpha(x,0) = \lambda$ (see Stage $0$ of the construction on page~96
of \cite{Cholak:94*1}), and if $x$ enters $A$ at stage $s$ then
$\alpha(x,s+1)\uparrow$ (see Step~1 on page~97), then, for almost all
$x$, there is a least stage $s_{\beta}$ such that either
$\alpha(x,s_\beta)\uparrow$ or $\beta \subseteq \alpha(x,s_\beta)$.
Let $R = \{ x | x \in A_{s_\beta}\}$.  $R$ is a computable subset of
$A$.  Assume $x \in \overline{R}$ enters $\Vhat^\dagger_{n}
=\Vhat^\dagger_{g(e)} = \Vhat^\dagger_\beta$ at stage $s$.  Let $s' =
\max\{s,s_\beta\}$.  By Equation~\eqref{eq:1} and the definition of
$R$, $x \not\in A_{s'}$ and hence $\beta \subseteq \alpha(x,s')$.
Then, by the last clause of $\mathcal{Q}_{\alpha}$ (on page 95), $x
\in \Vhat^{+}_{\beta,s'} = \Vhat^{+}_{g(e),s'} = \Vhat^{+}_{e,s'}$.
By  $\mathcal{Q}_{\alpha}$, $\Vhat^{+}_{e} \subseteq
\Vhat^\dagger_{n}$. Hence $\Vhat^{+}_{e} \doteq_\R \Vhat^\dagger_n$.
The proof that $\Uhat^{+}_{e} \doteq_\R \Uhat^\dagger_n$ is similar.
$\hfill \Box$

\subsection{Meeting the hypothesis of the Modified Translation
Theorem}\label{sec:meet-hypoth-modif}

By the hypothesis of Theorem~\ref{extension} and Definition~\ref{extiso},
we can assume that there are an effective listing of splits of $A$,
$\{S_i\}_{i \leq \omega}$, and a $\Delta^0_3$ set $B$ such that
$\{S_i\}_{i \in B}$ generates $\B$ and $\{\Shat_i\}_{i \leq \omega}$
is a similar listing of splits of $\Ahat$ for $\widehat{\mathcal{B}}$,
$\Bhat$, and $\Ahat$ such that $\Theta(S_i) = \Shat_{\Theta(i)}$ and
$\Theta^{-1}(\Shat_i) = \Shat_{\Theta^{-1}(i)}$ is an extendible
isomorphism between $\B$ and $\widehat{\B}$.

By Lemmas~\ref{trivialiso} and \ref{joinextiso}, we can assume that
the split $\emptyset$ and $A$ appears as some $S_i$ and $\breve{S}_i$
for some $i \in B$.  Since $\{S_i\}_{i \leq \omega}$ is effective we
can assume for all $i$, $S_{2i+1} = \breve{S}_{2i}$ and that $2i \in
B$ iff $2i+1 \in B$.  Similarly for $\{\Shat_i\}_{i \leq \omega}$ and
$\widehat{B}$.  Without loss, we can assume that
$S_{\Theta^{-1}(2e+1)} = \breve{S}_{\Theta^{-1}(2e)}$ and
$\Shat_{\Theta(2e+1)} = \breve{\Shat}_{\Theta(2e)}$.  Since
$\{S_i\}_{i \leq \omega}$ and $\{\Shat_i\}_{i \leq \omega}$ are
effective listings of splits, $\Theta$ remains $\Delta^0_3$.  By
Lemma~\ref{entryBA}, we will also assume that for all $i$, $A\searrow
S_i = \emptyset$, for some fixed enumeration of $\{A\}_{s \leq
  \omega}$. Dually for $\{\Shat_i\}_{i \leq \omega}$ and $\Ahat$.

Furthermore, since at this point we no longer need an effective
enumeration of splits, if $2i \not\in B$, let $S_{2i} = \emptyset$,
$\Shat_{\Theta(2i)} = \emptyset$, $S_{2i+1} = A$ (with the enumeration
$\{A_{s+1}\}_{s \in \omega}$ so $A^\dagger \searrow S_{2i+1}
=\emptyset$) and $\Shat_{\Theta(2i+1)} = \Ahat$ (with the enumeration
$\{\Ahat_{s+1}\}_{s \in \omega}$ so $\Ahat^\dagger \searrow
\Shat_{\Theta(2i+1)} =\emptyset$) and dually for $\{\Shat_i\}_{i \leq
  \omega}$ and $\widehat{B}$.

We want to, using an oracle for $\mathbf{0''}$, inductively construct
an enumeration of the c.e.\ sets $\{U^\dagger_{n}\}_{n< \omega}$,
$\{\Vhat^\dagger_{n}\}_{n < \omega}$, $\{\Uhat^\dagger_{n}\}_{n<
  \omega}$, and $\{V^\dagger_n\}_{n < \nu}$ which meets the two
hypotheses of Theorem~\ref{mtt}.  Let $\mathcal{N}_e$ be the set of
$(2e+1)$-states $\nu$ such that $D^{A}_\nu$ is infinite and
$D^{\Ahat}_\nu$ is infinite, where $D^{A}_\nu$ is measured w.r.t.\ 
$\{S_{n,s}\}_{n\leq 2e+1}$ and $\{S_{\Theta^{-1}(n),s}\}_{i\leq2e+1}$
and $D^{\Ahat}_\nu$ is measured w.r.t.\ 
$\{\Shat_{\Theta(n),s}\}_{n\leq 2e+1}$ and
$\{\Shat_{n,s}\}_{i\leq2e+1}$, for all $s < \omega$.  Determining
$\mathcal{N}_e$ is the only place $\mathbf{0''}$ is used.

Let $x \in A_{s+1}-A_{s}$.  Let $\nu = \nu(2e+1,x,s)$ (as measured
above).  If $\nu \in \mathcal{N}_e$ then let $x \in U^\dagger_{2e,s}$
iff $x \in S_{2e,s}$, $x \in U^\dagger_{2e+1,s}$ iff $x \in
S_{2e+1,s}$, $x \in \Vhat^\dagger_{2e,s}$ iff $x \in
S_{\Theta^{-1}(2e),s}$, and $x \in \Vhat^\dagger_{2e+1,s}$ iff $x \in
{S}_{\Theta^{-1}(2e+1),s}$. We act dually if $\xhat \in
\Ahat_{s+1}-\Ahat_{s}$.  For all $s$, let $A^\dagger_{s} = A_s$ and
$\Ahat^\dagger_{s} = \Ahat_s$.
 
Since only finitely much information, mainly $\mathcal{N}_e$, is used in
the above construction of the sets $U^\dagger_{2e}$,
$U^\dagger_{2e+1}$, $\Vhat^\dagger_{2e}$, $\Vhat^\dagger_{2e+1}$,
$\Uhat^\dagger_{2e}$, $\Uhat^\dagger_{2e+1}$, $V^\dagger_{2e}$, and
$V^\dagger_{2e+1}$, these sets are computably enumerable.  Hence
$\{U^\dagger_{n,s}\}_{n,s < \omega}$, $\{\Vhat^\dagger_{n,s}\}_{n,s <
  \omega}$, $\{\Uhat^\dagger_{n,s}\}_{n,s < \omega}$, and
$V^\dagger_{n,s}\}_{n,s < \omega}$, is a $\mathbf{0''}$-enumeration of
$\{U^\dagger_{n}\}_{n < \omega}$, $\{\Vhat^\dagger_{n}\}_{n <
  \omega}$, $\{\Uhat^\dagger_{n}\}_{n < \omega}$, and
$\{V^\dagger_{n}\}_{n < \omega}$ satisfying Condition~\eqref{eq:1}.
By induction on $e$, we can easily show that for all $(2e+1)$-states
$\nu$, $\nu \in \mathcal{N}_e$ iff $D^{A^\dagger}_\nu$ is infinite iff
$D^{\Ahat^\dagger}_{\nu}$ is infinite, where $D^{A^\dagger}_\nu$ and
$D^{\Ahat^\dagger}_{\nu}$ are measured as in Theorem~\ref{mtt}.
Therefore Condition~\eqref{eq:2} is satisfied.

\begin{lemma}\label{sec:computable}
  For all $e$, $U^\dagger_{2e} \doteq_\R S_{2e}$, $U^\dagger_{2e+1}
  \doteq_\R S_{2e+1}$, $\Vhat^\dagger_{2e} \doteq_\R
  S_{\Theta^{-1}(2e)}$, and $\Vhat^\dagger_{2e+1} \doteq_\R
  {S}_{\Theta^{-1}(2e+1)}$. For all $e$, $\Uhat^\dagger_{2e} \doteq_\R
  \Shat_{\Theta(2e)}$, $\Uhat^\dagger_{2e+1} \doteq_\R
  \Shat_{\Theta(2e+1)}$, $V^\dagger_{2e} \doteq_\R \Shat_{2e}$, and
  $V^\dagger_{2e+1} \doteq_\R \Shat_{2e+1}$.
\end{lemma}

\begin{proof}
  Since $\Theta$ is an isomorphism between $\B$ and $\widehat{\B}$,
  for each $(2e+1)$-state $\nu$, $\{x: \nu(2e+1,x) = \nu\}$ is
  noncomputable iff $\{\xhat : \widehat{\nu}(2e+1,\xhat) = \nu\}$ is
  noncomputable, where $\nu(2e+1,x)$ is measured w.r.t.\ 
  $\{S_i\}_{i\leq2e+1}$ and $\{S_{\Theta^{-1}(i)}\}_{i \leq 2e+1}$ and
  $\nu(2e+1,\xhat)$ is measured w.r.t.\ 
  $\{\Shat_{\Theta(i)}\}_{i\leq2e+1}$ and $\{\Shat_i\}_{i \leq 2e+1}$.
  
  By our carefully chosen enumerations of splits of $A$,
  $\{S_{i,s}\}_{i,s \leq \omega}$, the set $\{x: \nu = \nu(2e+1,x,s)
  \wedge x \in A_{s+1}-A_s\}$ is noncomputable iff $\{x: \nu=
  \nu(2e+1,x)\}$ is noncomputable, where $\nu(2e+1,x,s)$ is measured
  as above.  Dually for $\Ahat$.
  
  Let $\mathcal{A}_e$ be the set of all $(2e+1)$-states $\nu$.  For
  $\nu \in \mathcal{A}_e$, let $S_{2e,\nu}$ be the set $\{ x : \nu =
  \nu(2e+1,x,s) \wedge x \in A_{s+1}-A_s \wedge x \in S_{2e,s}\}$. If
  $\nu \not\in \mathcal{N}_e$ then $S_{2e,\nu}$ is computable.
  $S_{2e} = \bigsqcup_{\nu \in \mathcal{A}_e} S_{2e,\nu}$. By the
  above construction, $U^\dagger_{2e}= \bigsqcup_{\nu \in
    \mathcal{N}_e} S_{2e,\nu}$.  Hence $U^\dagger_{2e} \doteq_\R
  S_{2e}$.  We can argue similarly for the remaining sets.
\end{proof}

Since $\Theta$ is an isomorphism between substructures of $\S_\R(A)$
and $\S_\R(\Ahat)$, $A$ is noncomputable iff $\Ahat$ is
noncomputable.  As we noted shortly after the statement of
Theorem~\ref{extension}, Theorem~\ref{extension} holds when $A$ and
$\Ahat$ are computable.

\subsection{Constructing the isomorphism $\Phi$}

In the above section we built a $\mathbf{0''}$-enumeration meeting the
hypothesis of Theorem~\ref{mtt} and satisfying
Lemma~\ref{sec:computable}.  Now apply Theorem~\ref{mtt} to this
enumeration.  Conditions~\eqref{eq:4}, \eqref{eq:7}, and \eqref{eq:8}
of Theorem~\ref{mtt} are the three conditions in the hypothesis of
Soare's original Extension Theorem (see \citet{Soare:87}
Theorem~XV.4.5).  Now apply Soare's original Extension Theorem to the
enumeration given to us by Theorem~\ref{mtt}.  This gives us the c.e.\ 
sets $\{U_{n}\}_{n < \omega}$, $\{\Vhat_{n}\}_{n < \omega}$,
$\{\Uhat_{n}\}_{n < \omega}$, and $\{V_{n}\}_{n < \omega}$. The
Extension Theorem only adds elements to $\Vhat^{+}_n$ to get
$\Vhat_n$ and similarly for $\Uhat_n$.  $\Phi(U_n) =\Uhat_n$ and
$\Phi^{-1}(V_n) = \Vhat_n$ is an isomorphism between $\E^*(A)$ and
$\E^*(\Ahat)$ (see \citet{Soare:87}~Section XV.4 for details).

By Lemma~\ref{sec:computable}, for all $n$, $U^\dagger_n \doteq_\R
S_n$ and $\Uhat^\dagger_n \doteq_\R \Shat_{\Theta(n)}$.  By
\eqref{eq:5} of Theorem~\ref{mtt}, $U_{e_n} =^* U^\dagger_n$ and
$\Uhat^{+}_{e_n} \doteq_\R \Uhat^\dagger_{n}$. Therefore for all $n$,
$U_{e_n} \doteq_\R S_n$, $\Uhat^{+}_{e_n} \doteq_\R
\Shat_{\Theta(n)}$.  Since $\Theta$ is an isomorphism,
$\Theta(U_{e_n}) =_\R \Theta(S_n)$.

By our careful choice of $\{S_i\}_{i < \omega}$ and our modification
of $\Theta$ in Section~\ref{sec:meet-hypoth-modif} we have that for
all $n$, $S_{2n} \sqcup S_{2n+1} = A$ and $\Shat_{\Theta(2n)} \sqcup
\Shat_{\Theta(2n+1)} = \Ahat$.  Hence for all $n$, $U_{e_{2n}} \sqcup
U_{e_{2n+1}} \sqcup R_n = A$ and $\Uhat^{+}_{e_{2n}} \sqcup
\Uhat^{+}_{e_{2n+1}} \sqcup \widehat{R}_n = \Ahat$, for some
computable sets $R_n$ and $\widehat{R}_n$.

Since $\Phi$ is an isomorphism between $\E^*(A)$ and $\E^*(\Ahat)$ and
the sets $S_{2n}$ and $U_{e_{2n+1}}$ are disjoint, $\Uhat_{e_{2n}} -
\Uhat^{+}_{e_{2n}} \subseteq^* \widehat{R}_n$ and $\Phi(S_{2n})
- \Uhat_{e_{2n}} \subseteq^* \Rhat_n$.  Therefore $\Phi(U_{e_{2n}})
=^* \Uhat_{e_{2n}} =_\R \Shat_{\Theta(2n)}$ and $\Phi(S_{2n}) =_\R
\Phi(U_{e_{2n}})$.  So $\Phi(S_{{2n}}) =_\R \Theta(S_{{2n}})$.
We argue similarly to show $\Phi(S_{2n+1}) =_\R \Theta(S_{2n+1})$ and
$\Phi^{-1}(\Shat_{n}) =_\R \Theta^{-1}(\Shat_{n})$.  $\hfill \Box$

\section{Extensions to automorphisms}\label{ext2} %

Our goal to find an algebraic extension theorem which allows us to
find an automorphism $\Lambda$ of $\E$ taking $A$ to $\Ahat$ if and
when possible.  Clearly we will have to add some extra hypotheses to
Theorem~\ref{extension} about the outside of $A$ and $\Ahat$.

Recall that $\L^*(A)$ is the structure $(\{W_e \cup A : e \in
\omega\}, \subseteq)$ modulo the finite sets.  A substructure $\L$
of $\L^*(A)$ is a subcollection of the sets $(\{W_e \cup A : e \in
\omega\}, \subseteq)$ modulo the finite sets.  An isomorphism between
$\L^*(A)$ and $\L^*(\Ahat)$ is a one-to-one, onto (both of these items
are in terms of $*$-equivalence classes) function $\Xi$ from $\{W_e :
e \in \omega\}$ to $\{\What_e : e \in \omega\}$ such that $W_e \cup A
\subseteq^* W_i \cup A$ iff $\Xi(W_e \cup \Ahat) \subseteq^* \Xi(W_i
\cup \Ahat)$.  Note that $\Xi$ is applied to $W \cup A$.

Assume that $\L^*(A)$ and $\L^*(\Ahat)$ are isomorphic via $\Psi$ and
that $\B$ and $\widehat{\B}$ are $\Delta^0_3$ isomorphic via $\Theta$.
We wish to use the isomorphism $\Phi$ from Theorem~\ref{extension} to
extend this pair of isomorphisms into an automorphism $\Lambda$ of
$\E$ such that $\Lambda(A) = \Ahat$.

Notice that $W = (W - A) \sqcup (W \cap A)$.  It would be nice to
define $\Lambda(W) = (\Psi(W \cup A) - \Ahat) \sqcup \Phi(W \cap A)$.
Clearly this is order preserving.  But why is $(\Psi(W \cup A) -
\Ahat) \sqcup \Phi(W \cap A)$ a \ce set? To answer that we must
explore more carefully the complex relation between $\L^*(A)$ and
$\B$.

\begin{definition}
  $S$ \emph{supports} $X$ iff $S \subseteq X$ and $(X-A) \sqcup S$ is
  a \ce set.
\end{definition}

\begin{lemma}\label{supportsisSigma03}
   Whether $S$ supports $X$ is $\Sigma^0_3$.
\end{lemma}

\begin{proof}
  $S$ supports $X$ iff there exists an $e$ where $W_e = (X-A) \sqcup
  S$ and $S \subseteq X$.
\end{proof}

\begin{lemma}
  $W \searrow A$ supports $W$.
\end{lemma}

\begin{proof}
  $W = (W -A) \sqcup (W \searrow A) \sqcup (A \searrow W)$ and $(W -A)
  \sqcup (W \searrow A)$ is the \ce set $W \backslash A$.
\end{proof}

\begin{definition}
  An extendible subalgebra $\B$ \emph{supports} $\L$ if for all $W \in
  \L$ there an $i \in B$ such that $S_i$ supports $W$.%
\end{definition}

\begin{lemma} \label{supports}
  $\E_A$ supports $\L^*(A)$.
\end{lemma}

\begin{lemma} \label{psupports}
  If $S$ supports $X$ and $T$ is a split of $A$ such that $T \subseteq
  S$ and $S =_{\R(A)} T$ then $T$ supports $X$.
\end{lemma}

\begin{proof}
  $(X-A) \sqcup S$ is a \ce set.  If $S-T$ is a computable set $R$ then
  $(X-A) \sqcup T = (((X-A) \sqcup S) \cap \overline{R})$ is a \ce
  set.
\end{proof}

\begin{definition}
  Assume that
  \begin{itemize}
  \item $\L^*(A)$ and $\L^*(\Ahat)$ are isomorphic via $\Psi$,
  \item $\B$ and $\widehat{\B}$ are isomorphic via $\Theta$,
  \item $\B$ {supports} $\L$, and
  \item $\widehat{\B}$ {supports} $\widehat{\L}$.
  \end{itemize}
  Then the isomorphisms $\Psi$ and $\Theta$ \emph{preserve} the
  supports of $\L$ and $\widehat{\L}$ if
   \begin{itemize}
   \item for $W^* \in \L$, there is an $i \in B$ such that $S_i$
     supports $W$ and $(\Psi(W \cup A) - \Ahat) \sqcup \Theta(S_i)$ is
     a \ce set, and
   \item for all $\What^* \in \widehat\L$, there is an $i \in
     \widehat{B}$ such that $\Shat_i$ supports $\What$ and
     $(\Psi^{-1}(\What \cup \Ahat) - A) \sqcup \Theta^{-1}(\Shat_i)$
     is a \ce set.
  \end{itemize}
  For shorthand we just say isomorphisms $\Psi$ and $\Theta$
  \emph{preserve supports}.
\end{definition}

If $S_i$ supports $W$ then $S_i \subseteq W$.  But if isomorphisms
$\Psi$ and $\Theta$ preserve supports, then, while $(\Psi(W \cup A) -
\Ahat) \sqcup \Theta(S_i)$ is a \ce set, we do not require that
$\Theta(S_i)$ be contained in $\Psi(W)$.  Hence $\Theta(S_i)$ might not
be a support of $\Psi(W)$.

\begin{theorem}\label{extension3}  Assume that
  \begin{enumerate}
  \item $\L^*(A)$ and $\L^*(\Ahat)$ are isomorphic via $\Psi$,
  \item $\B$ and $\widehat{\B}$ are extendible algebras which are
    extendibly $\Delta^0_3$ isomorphic via $\Theta$,
  \item $\B$ supports $\L^*(A)$,
  \item $\widehat{\B}$ supports $\L^*(\Ahat)$,
  \item $\Psi$ and $\Theta$ preserves supports,
  \item \label{needed} $\Phi$ is an isomorphism between $\E^*(A)$ and
    $\E^*(\Ahat)$ such that if $i \in B$ then $\Theta(S_i) =_{\R}
    \Phi(S_i)$ and if $i \in \widehat{B}$ then $\Theta^{-1}(\Shat_i)
    =_{\R} \Phi^{-1}(\Shat_i)$.
  \end{enumerate}
  Then $\Lambda(W) = (\Psi(W \cup A) - \Ahat) \sqcup \Phi(W \cap A)$
  is an automorphism of $\E$ taking $A$ to $\Ahat$.
\end{theorem}

\begin{proof}
  It is enough to show that $(\Psi(W \cup A) - \Ahat) \sqcup \Phi(W
  \cap A)$ is a \ce set.  First note that $W \cap A = S_i \sqcup
  (\breve{S_i} \cap W)$, where $S_i$ supports $W$ and $i \in B$.
  Since $\Phi$ is an isomorphism between $\E^*(A)$ and $\E^*(\Ahat)$,
  $\Phi(W \cap A) = \Phi(S_i) \sqcup \Phi(\breve{S_i} \cap W)$.  Since
  $\Psi$ and $\Theta$ preserve supports, for some support $S_i$ of
  $W$, $(\Psi(W \cup A) -\Ahat) \sqcup \Theta(S_i)$ is a \ce set.
  Since $\Theta(S_i) =_{\R} \Phi(S_i)$, $(\Psi(W \cup A) - \Ahat)
  \sqcup \Phi(S_i)$ is a \ce set. Hence $(\Psi(W \cup A) - \Ahat)
  \sqcup \Phi(S_i) \sqcup \Phi(\breve{S}_{i} \cap W)$ is a \ce set.
  Similarly we can show $\Lambda^{-1}(\What)$ is a \ce set.
\end{proof}

\begin{theorem}\label{extension2}  Assume that
  \begin{enumerate}
  \item $\L^*(A)$ and $\L^*(\Ahat)$ are isomorphic via $\Psi$,
  \item $\B$ and $\widehat{\B}$ are extendible algebras which are
    extendibly $\Delta^0_3$ isomorphic via $\Theta$,
  \item $\B$ supports $\L^*(A)$,
  \item $\widehat{\B}$ supports $\L^*(\Ahat)$,
  \item $\Psi$ and $\Theta$ preserve supports.
 \end{enumerate}
 Then there is an automorphism $\Lambda$ of $\E$ such that $\Lambda(A)
 = \Ahat$, $\Lambda \restriction \L^*(A) = \Psi$, and $\Lambda
 \restriction \E^*(A)$ is $\Delta^0_3$.
\end{theorem}

\begin{proof}
  Apply Theorem~\ref{extension} to get $\Phi$ as required by
  Theorem~\ref{extension3}.\ref{needed}.  $\Phi$ is $\Delta^0_3$.
  Apply Theorem~\ref{extension3} to get $\Lambda$.
\end{proof}

The way we put together the automorphism in Theorem~\ref{extension2}
is very similar to the way in which Herrmann showed that the Herrmann sets
(along with the hemimaximal sets and other such orbits) form an orbit
(see \citeasnoun[Sections 5 and 6]{Cholak.Downey.ea:01}).  Both
methods are algebraic or ``static''.

In Section~\ref{mainsec}, we will show that Theorem~\ref{extension2}
can be improved to be an ``if and only if' statement (see
Theorem~\ref{orbits}).

\section{Preserving the computable subsets}\label{computable}

\begin{definition}
  A map $\Xi$ from a substructure of $\mathcal{G} \subseteq \E(A)$ to
  $\widehat{\mathcal{G}} \subseteq \E(\widehat{A})$ \emph{preserves
    the computable subsets} if $R \in \R(A) \cap \mathcal{G}$ iff
  $\Xi(R) \in \R(\Ahat) \cap \widehat{\mathcal{G}}$.
\end{definition}

There is no guarantee that any of the maps we have been considering
preserves the computable subsets; this includes $\Theta$. And the same
can be said about Soare's original Extension Theorem (see
\citeasnoun[XV.4.5]{Soare:87}) (applied by itself).  To see this: 
If $X \in \R(A)$ and $\Theta$ is an isomorphism $\Theta$ between
$\E^*(A)$ and $\E^*(\Ahat)$, then there is a $Y$ such that $X \sqcup Y
= A$ and $\Theta(X) \sqcup \Theta(Y) = \Ahat$ but there may not be a
$Z$ such that $\Theta(X) \sqcup Z = \widehat{\omega}$.  Of course,
there is such a $Z$ if $\Ahat$ is computable (and dually if $A$ is
computable).

It might be useful to consider the following example: If $A$ and
$\Ahat$ are infinite then there is an effective isomorphism $\Psi$
between $\E^*(A)$ and $\E^*(\Ahat)$ (let $f$ be an effective map from
$A$ to $\Ahat$ and let $\Psi(W) = f(W)$).  If $A$ is computable but
$\Ahat$ is not then $\Psi$ cannot preserve the computable subsets.

From this point on we will always consider $A$ and $\Ahat$ to be
noncomputable.  We will point out that it is known that there is an
isomorphism between $\E^*(A)$ and $\E^*(\Ahat)$ which preserves the
computable subsets (see Theorem~\ref{apcomputable}). The goal of this
section is to provide another proof of fact using our methods.

\begin{definition}
  $\C(A)$ is the set of $W_e$ such that either $\overline{A} \subseteq
  W_e$ or $ W_e \subseteq^* A$.
\end{definition}

\begin{theorem}[Soare's Automorphism Theorem \cite{Soare:74}] 
   \label{apcomputable}
   Let $A$ and $\Ahat$ be two noncomputable \ce sets.
   \begin{enumerate}
   \item Then there is a $\Delta^0_3$ isomorphism $\Lambda$ between $\E(A)
     \cup \C(A)$ and $\E(\Ahat) \cup \C(\Ahat)$. Furthermore a
     $\Delta^0_3$-index for $\Lambda$ can be found uniformly from
     indexes for $A$ and $\Ahat$.
   \item In addition, $\Lambda$ preserves the computable subsets of $A$.
   \end{enumerate}
\end{theorem}

\citet{Soare:74} explicitly stated Theorem~\ref{apcomputable}.1.
Soare's result that maximal sets are automorphic follows since $A$ is
maximal iff $\C(A) = \E^*$. 

Theorem~\ref{apcomputable}.2 was observed, in unpublished work, by
Herrmann.  Assume that $R$ is a computable subset of $A$. Herrmann's
observation was that $\overline{R} \in \C(A)$ and hence $\Lambda(R)
\sqcup \Lambda(\overline{R}) =^* \widehat{\omega}$ and therefore
$\Lambda$ maps $R$ to a computable subset of $\Ahat$.  This
observation of Herrmann was never published and is one of the key
facts he used in showing that the Herrmann sets form an orbit; see
\citet{Cholak.Downey.ea:01}.

\subsection{Another proof of Theorem~\ref{apcomputable}}

We would like to show Theorem~\ref{apcomputable} using the methods of
this paper. 

First note that an isomorphism $\Lambda$ between $\E^*(A)$ and
$\E^*(\Ahat)$ preserving the computable subsets induces an isomorphism
$\Lambda'$ between $\E^*(A) \cup \C(A)$ and $\E^*(\Ahat) \cup
\C(\Ahat)$ taking $A$ to $\Ahat$.  If $\overline{A} \subseteq W$ then
$A \cup W = \omega$ and there is a computable set $R \subseteq A$ ($R
= A \backslash W$) such that $\overline{R} \subseteq W$ which implies
$W = \overline{R} \sqcup (W \cap R)$.  So for $W \in \mathcal{C}(A)$,
let $\Lambda'(W)$ be $\overline{\Lambda(R)} \sqcup (\Lambda(W \cap
R))$.

We would like to prove a theorem along the lines of
Theorem~\ref{extension3}.

\begin{theorem}\label{extension4}
  Assume that
  \begin{enumerate} 
  \item $\B$ and $\widehat{\B}$ are extendible algebras which are
    $\Delta^0_3$ extendibly isomorphic via $\Theta$;
  \item for all $R \in \R(A)$, there is an $i \in B$ such that $S_i$
    is computable and $R \subseteq S_i$;
  \item for all $\widehat{R} \in \R(\Ahat)$, there is an $i \in
    \widehat{B}$ such that $\Shat_i$ is computable and $\widehat{R}
    \subseteq \Shat_i$;
  \item for all $i \in {B}$, $\Theta(S_i)$ is computable iff $S_i$ is
    computable and for all $i \in \widehat{B},$ $\Theta^{-1}(\Shat_i)$
    is computable iff $\Shat_i$ is computable.
 \end{enumerate}
 Then there is a $\Lambda$ such that $\Lambda$ is a $\Delta^0_3$
 isomorphism between $\E^*(A)$ and $\E^*(\Ahat)$ which preserves the
 computable subsets, for all $i \in B$, $\Lambda(S_i) =_{\R}
 \Theta(S_i)$, and if $i \in \widehat{B}$, then $\Lambda^{-1}(\Shat_i)
 =_\R \Theta^{-1}(\Shat_i)$.
\end{theorem}

\begin{proof}
  First apply Theorem~\ref{extension} to get $\Phi$.  We will show
  that $\Phi$ is the desired isomorphism $\Lambda$.  It is enough to
  show $\Phi$ preserves the computable subsets.
  
  Let $R \in \R(A)$.  There is an $i$ such that $S_i$ is computable
  and $R \subseteq S_i$.  $\Theta(S_i)$ is computable.  By
  Theorem~\ref{extension}, $\Phi(S_i) \equiv_{\R} \Theta(S_i)$.  Hence
  $\Phi(S_i)$ is computable.  Therefore, since the set
  $A-R$ is c.e., the set $\overline{\Phi(R)} =^* \overline{\Phi(S_i)}
  \sqcup \Phi(S_i \cap (A -R))$ is \ce and $\Phi(R)$ is computable.
  The other direction is similar.
\end{proof}

It is actually reasonably easy to meet the hypothesis of the above
theorem; it is enough that $A$ and $\Ahat$ both be noncomputable.

\begin{theorem}\label{meeting}
  Let $A$ and $\Ahat$ be two noncomputable \ce sets.  Then there are
  $\B$ and $\widehat{\B}$ such that
   \begin{enumerate} 
   \item $\B$ and $\widehat{\B}$ are extendible algebras which are
     $\Delta^0_3$ extendibly isomorphic via $\Theta$;
   \item for all $R \in \R(A)$, there is an $i \in B$ such that $S_i$
     is computable and $R = S_i$;
   \item for all $\widehat{R} \in \R(\Ahat)$, there is an $i \in
     \widehat{B}$ such that $\Shat_i$ is computable and $\widehat{R} =
     \Shat_i$;
   \item for all $i \in {B}$, $\Theta(S_i)$ is computable iff $S_i$ is
     computable and for all $i \in \widehat{B},$
     $\Theta^{-1}(\Shat_i)$ is computable iff $\Shat_i$ is computable.
 \end{enumerate} 
\end{theorem}

\begin{proof}
  Apply Theorem~\ref{sec:extend-subalg} and its dual to get $\B$ and
  $\widehat{\B}$.  Now both $B$ and $\widehat{B}$ are infinite and
  $\Delta^0_3$.  We will inductively define $\theta$.  If $i +1 \in
  B$, let $\theta({i+1})$ be the least element of $\widehat{B}$ which
  is not yet in the range of $\theta$.  Otherwise $\theta(i+1)$ is
  undefined.  Let $\Theta(S_i) = \Shat_{\theta(i)}$.  Similarly for
  $\Theta^{-1}$.  Clearly $\Theta$ is $\Delta^0_3$.
  
  Since everything in $B$ and $\widehat{B}$ are computable splits of
  $A$, $\B$ and $\widehat{\B}$ are classically isomorphic to the
  trivial Boolean algebra.  Therefore $\Theta$ induces an isomorphism
  between $\B$ and $\widehat{\B}$.  Hence $\Theta$ is clearly the
  desired extendible isomorphism.
\end{proof}

By combining Theorems~\ref{extension4} and \ref{meeting} we get
another proof of Theorem~\ref{apcomputable}.

\subsection{Some examples of the use of 
  Theorem~\ref{apcomputable}}\label{examples}

\subsubsection{The hemimaximal sets}  \label{hemimaximal}

We include this example as it has not appeared previously in print in
this form and it hints of things to come in later sections.  Assume
$A_1\sqcup A_2 = A$ where the $A_i$s are not computable.  Dually for
$\Ahat$. Assume that $\Theta_i$ is an isomorphism from $\E^*(A_i)$ to
$\E^*(\Ahat_i)$ that preserves the computable subsets (from
Theorem~\ref{apcomputable}).

As with the maximal sets, it is enough to define an isomorphism
$\Lambda$ between $\E^*(A)$ and $\E^*(\Ahat)$ preserving the
computable subsets. If $X \subseteq^* A$ then let $\Lambda(X) =
\Theta_1(X \cap A_1) \sqcup \Theta_2(X \cap A_2)$.  Let ${R} \in
{\R}(A)$. Then ${R} \cap A_i$ is computable.  So $\Theta_i({R} \cap
A_i)$ is computable.  Hence $\Theta_1({R} \cap A_1) \sqcup
\Theta_2({R} \cap A_2)$ is computable.  The complexity of the
resulting automorphism is $\Delta^0_3$.

Downey and Stob's proof used the fact that if $W \cup A = \omega$ then
$W \searrow A_i$ is infinite: a very dynamic property.  Our proof only
relies on algebraic facts.

\subsubsection{The atomless Boolean Algebra $\S_{\R}(A)$}

As we know, all atomless Boolean Algebras are isomorphic but with
$\S_{\R}(A)$ something stronger is true.

\begin{theorem}[Nies, see \citet{MR2004f:03077}] %
  \label{nies}
  If $A$ and $\widehat{A}$ are noncomputable, then $\S_\R(A)$ and
  $\S_\R(\widehat{A})$ are $\Delta^0_3$ isomorphic.
\end{theorem}

\begin{proof}
  The isomorphism $\Lambda$, from Theorem~\ref{apcomputable}, is an
  isomorphism between $\E^*(A)$ and $\E^*(\Ahat)$ preserving the
  computable sets.  Hence $\Lambda$ induces an isomorphism between
  $\S_\R(A)$ and $\S_\R(\widehat{A})$.
\end{proof}

\subsection{Extendible Algebras of Computable Sets}

This section was added after the rest of the paper was completed.  As
we mentioned in the Introduction (third to last paragraph) and last
sentence, this paper has a sequel.  The goal of this section is to
provide a clear, clean interface between the two papers.  In
particular, we will proof a theorem, Theorem~\ref{interface}, which we
hope we can use as a black box in the sequel.  

Theorem~\ref{interface} is an improved version of
Theorem~\ref{apcomputable}.  In Theorem~~\ref{apcomputable} the
computable sets are preserved.  In Theorem~\ref{interface} the
computable sets are preserved plus an external isomorphism determines
where some of the computable sets are mapped.

\begin{definition}
  An extendible algebra $\B$ of $\S_\R(\omega)$ is called an
  \emph{extendible algebra of computable sets} as the splits of
  $\omega$ are the computable sets.
\end{definition}

\begin{lemma}
  If $\B = \{R_i : i \in B\} $ is an extendible algebra of computable
  sets then $\B_A = \{ R_i \cap A: i \in B\}$ is an extendible
  algebra of $\S_\R(A)$.
\end{lemma}

\begin{proof}
  $\{\tilde{R_i} \cap A: i \in \omega\}$ witnesses that $\{{R_i}
  \cap A: i \in \omega\}$ is an effective listing of splits of $A$.
\end{proof}

\begin{lemma}
  Assume that $\B$ and $\widehat{\B}$ are extendible subalgebras of
  computable sets which are extendibly isomorphic via $\Pi$.  $\Pi_A(R
  \cap A) = \Pi(R) \cap \Ahat$ is an extendible isomorphism between
  $\B_A$ and $\B_{\Ahat}$.
\end{lemma}

\begin{theorem}\label{interface}
  Let $\B$ be a extendible algebra of computable sets and similarly
  for $\widehat{\B}$.  Assume the two are extendibly isomorphic via
  $\Pi$.  Then there is a $\Phi$ such that $\Phi$ is a $\Delta^0_3$
  isomorphism between $\E^*(A)$ and $\E^*(\widehat{A})$, $\Phi$ maps
  computable subsets to computable subsets, and, for all $R \in \B$,
  $(\Pi(R) - \widehat{A}) \sqcup \Phi(R \cap A)$ is computable\ (and
  dually).
\end{theorem} \medskip

\begin{proof}
  Apply Lemmas~\ref{joinext} and \ref{joinextiso} to $\B_A$,
  $\widehat{\B}_{\Ahat}$, $\Pi_A$, and the extendible algebras and
  extendible isomorphism from Theorem~\ref{meeting} to get
  $\tilde{\B}$, $\widehat{\tilde{\B}}$ and $\tilde{\Theta}$. Now apply
  Theorem~\ref{extension} to get $\Phi$.  By the proof of
  Theorem~\ref{extension4}, $\Phi$ preserves the computable sets.

  Since $\Pi$ is an isomorphism between extendible algebras of
  computable sets, $\Pi(R)$ is a computable set.  By
  Theorem~\ref{extension}, $\Theta(R\cap A) \triangle \Pi_A(R) = R_0$
  is a computable subset of $\Ahat$.  Since $\Theta(R \cap A)$ is a
  split of $\Ahat$, $\Theta(R \cap A) \cap R_0 = R_1$ is a computable
  subset of $\Ahat$.  Similarly, $\Pi_A(R) \cap R_0 = R_2$ is a
  computable subset of $\Ahat$. So $\Phi(R \cap A) = (\Pi_A(R) \sqcup
  R_1) \cap \overline{R}_2$. Hence
  \begin{equation*}
    \begin{split}
       (\Pi(R) \sqcup R_1) \cap \overline{R}_2  = & \big ( (\Pi(R) -
      \Ahat) \sqcup \Pi_A(R)  \sqcup R_1 \big ) \cap
      \overline{R_2} \\ = &(\Pi(R) - \widehat{A}) \sqcup \Phi(R \cap A).\\
    \end{split}
  \end{equation*} 
  So $(\Pi(R) - \widehat{A}) \sqcup \Phi(R \cap A)$ is computable as
  desired.  The dual is proved in a similar fashion.
\end{proof}

\section{Automorphisms back to automorphisms}\label{mainsec}

Assume that $A$ and $\Ahat$ are automorphic via $\Psi$.  Hence
$\L^*(A)$ and $\L^*(\Ahat)$ are isomorphic via $\Psi$.  Since $A$ and
$\Ahat$ are automorphic, the structures $\mathcal{S}_{\mathcal{R}}(A)$
and $\mathcal{S}_{\mathcal{R}}(\Ahat)$ are isomorphic structures
(since they are definable structures). In fact, from
\citeasnoun{MR2004f:03077}, we know much more is true.

\begin{theorem}[\textbf{The Restriction Theorem}; 
Theorem~1.2 of \citeasnoun{MR2004f:03077}]
   \label{oldsplits}
   If $A$ and $\Ahat$ are automorphic via $\Psi$ then the structures
   $\mathcal{S}_{\mathcal{R}}(A)$ and
   $\mathcal{S}_{\mathcal{R}}(\Ahat)$ are $\Delta^0_3$-isomorphic
   structures via an isomorphism $\Gamma$ induced by $\Psi$.
\end{theorem}

In other words there is an isomorphism $\Gamma$ between
$\mathcal{S}_{\mathcal{R}}(A)$ and $\mathcal{S}_{\mathcal{R}}(\Ahat)$
such that for all splits of $A$, $\Gamma(S) =_\R \Psi(S)$; for all
splits $\Shat$ of $\Ahat$, $\Gamma^{-1}(\Shat) =_\R \Psi^{-1}(\Shat)$;
and a $\Delta^0_3$-function $f$ such that for $W_e \in \S(A)$,
$W_{f(e)} =_\R \Gamma(W_e)$. (For more about this theorem we direct
the reader to \citeasnoun{MR2004f:03077}.)

\begin{theorem}\label{hidden}
  Assume $A$ and $\Ahat$ are automorphic via $\Psi$. Let $\tilde{\B}$
  be an extendible algebra (of $\S_\R(A)$).  Then there are extendible
  $\widehat{\B}$ (of $\S_\R(\Ahat)$) and $\Theta$ such that
  \begin{enumerate}
  \item $\widehat{\B}$ and $\tilde{\B}$ are extendibly
    $\Delta^0_3$-isomorphic via $\Theta$,
  \item if $i \in \tilde{B}$ and $S_i$ supports $W$ then $\Theta(S_i)$
    supports $\Psi(W)$.
  \end{enumerate}
\end{theorem}

The proof of this theorem appears in Section~\ref{neil2proof}.  We
should note that we must argue dynamically in this proof. We can
use this result to show the following theorem.

\begin{theorem}[\textbf{The Conversion Theorem}]\label{mainnew}
  If $A$ and $\Ahat$ are automorphic via $\Psi$ then they are
  automorphic via $\Lambda$ where $\Lambda \restriction \L^*(A) =
  \Psi$ and $\Lambda \restriction \E^*(A)$ is $\Delta^0_3$.
\end{theorem}

\begin{proof} 
  $\L^*(A)$ and $\L^*(\Ahat)$ are isomorphic via $\Psi$.  Recall from
  Lemma~\ref{entry}, $\mathcal{E}_A$ is the extendible algebra
  generated by the entry sets. Recall from Lemma~\ref{supports},
  $\mathcal{E}_A$ supports $\mathcal{L}^*(A)$.  Apply
  Theorem~\ref{hidden} to $\mathcal{E}_A$ to get
  $\widehat{\mathcal{E}}_A$ and $\Theta_A$ and dually to
  $\mathcal{E}_{\Ahat}$ to get $\widehat{\mathcal{E}}_{\Ahat}$ and
  $\Theta_{\Ahat}$. By Lemmas~\ref{joinext} and \ref{joinextiso}, $\B
  = \mathcal{E}_A \join \widehat{\mathcal{E}}_{\Ahat}$ and
  $\widehat{\B} = \widehat{\mathcal{E}}_{A} \join \mathcal{E}_{\Ahat}$
  are extendible algebras $\Delta^0_3$-isomorphic via $\Theta$. Since
  $\E_A$ supports $\L^*(A)$, $\B$ does too.  Similarly for
  $\widehat{\B}$ and $\L^*(\Ahat)$.  By the last property of
  Theorem~\ref{hidden}, isomorphisms $\Psi$ and $\Theta$ preserve
  supports.  Now apply Theorem~\ref{extension2}.
\end{proof} 

Also using Theorem~\ref{extension2} we can algebraically describe an
orbit of $A$.

\begin{theorem}\label{orbits}
  The \ce sets $A$ and $\Ahat$ are automorphic iff there are $\Psi$,
  $\B$, $\widehat{\B}$, and $\Theta$ such that
\begin{enumerate}
\item $\L^*(A)$ and $\L^*(\Ahat)$ are isomorphic via $\Psi$,
\item $\B$ and $\widehat{\B}$ are extendible algebras which are
  extendibly $\Delta^0_3$ isomorphic via $\Theta$,
\item $\B$ supports $\L^*(A)$,
\item $\widehat\B$ supports $\L^*(\Ahat)$,
\item the isomorphisms $\Psi$ and $\Theta$ preserve supports.
  \end{enumerate}
\end{theorem}

\subsection{Proof of Theorem~\ref{hidden}}\label{neil2proof}

To make life notationally easier we will prove the dual. So let
$\tilde{\B}$ be an extendible algebra of $\S_\R(\Ahat)$ and we will
build $\B$. 

By Theorem~\ref{oldsplits}, $\tilde{\B}$ and $\Gamma^{-1}(\tilde{\B})$
are $\Sigma^0_3$ algebras which are $\Delta^0_3$ isomorphic via
$\Gamma^{-1}$.  But $\Delta^0_3$ images and preimages of extendible
algebras need not be extendible. Hence we cannot let $\B =
\Gamma^{-1}(\tilde{\B})$.  We will construct $\B$ to be extendible and
extendibly isomorphic to $\tilde{\B}$ via $\Theta$ (and hence
isomorphic to $\Gamma^{-1}(\tilde{\B})$).  In fact we are going to
show something stronger; we will show $\E_A \join \B$ is isomorphic to
$\Gamma({\E}_A) \join \tilde{\B}$.

We are going to construct $\B$ and $\Theta$ via a standard tree
agreement.  We will construct a tree, $Tr$. At each node $\alpha$ of
the tree, we will construct the splits of $A$, $S_\alpha$ and
$\breve{S}_\alpha$.  We are going to build these splits as entry sets.
So for all $\alpha$, if $x$ enters $A$ at stage $s+1$ then  $x$
enters either $S_\alpha$ or $\breve{S}_\alpha$ at stage $s$.

The list $\{S_\alpha\}_{\alpha \in Tr}$ is an effective listing of
splits. $B = \{\alpha | \alpha \subset f \wedge |\alpha| \in
\tilde{B}\}$ is a $\Delta^0_3$ set.  So an extendible algebra, $\B$,
is created.

If $i \in \tilde{B}$ then let $\Theta(S_\alpha) = \tilde{S}_i$ and
$\Theta^{-1}(\tilde{S}_i) = S_\alpha$, where $\alpha \subset f$ and
$|\alpha| = i$. If we can show $\Theta$ induces an isomorphism between
$\B$ and $\tilde{\B}$ then $\Theta$ will be a $\Delta^0_3$-extendible
isomorphism between $\B$ and $\tilde{\B}$.  Hence without loss we can
assume that if $i \notin \tilde{B}$ then $\tilde{S}_i = \emptyset$ and
$\Gamma^{-1}(\tilde{S}_i) = \emptyset$.

For the rest of this proof we will use \emph{$e$-splits states} rather
than $e$-states.

\begin{definition}
  \begin{enumerate}
  \item For any $e$, if we are given a uniform enumeration of splits
    of $A$ $\{S_{i,s}\}_{i \leq e,s< \omega}$, $\{\breve{S}_{i,s}\}_{i
      \leq e,s< \omega}$, $\{T_{i,s}\}_{i \leq e, s < \omega}$, and
    $\{\breve{T}_{i,s}\}_{i \leq e, s < \omega}$ define the
    \emph{$e$-split state of $x$ at stage $s$}, $\nu^S(e,x,s)$, to be
    the full $2e$-state of $x$ w.r.t.\ $\{X_{i,s}\}_{i\leq 2e, s <
      \omega}$ and $\{Y_{i,s}\}_{i\leq 2e, s < \omega}$, where $
    X_{2i,s} = S_{i,s}$, $X_{2i+1,s} =\breve{S}_{i,s}$, $Y_{2i,s}=
    T_{i,s}$, and $ Y_{2i+1,s} = \breve{T}_{i,s} $.
  \item Let $\nu^S(\alpha,x,s) = \nu^S(e,x,s)$ where $|\alpha| = e$
    and $\nu^S(e,x,s)$ is measured w.r.t.\ $\{(W_i \searrow A)_s\}_{i
      \leq e,s<\omega}$, $\{(A \backslash W_i)_s\}_{i \leq
      e,s<\omega}$, $ \{S_{\beta,s}\}_{\beta \subseteq \alpha,s<
      \omega}$, and $\{\breve{S}_{\beta,s}\}_{\beta \subseteq
      \alpha,s< \omega}$.
  \item For any collection of splits of $A$, $\{S_{i}\}_{i \leq e}$
    and $\{T_{i}\}_{i \leq e}$, define the \emph{final $e$-split state
      of $x$} %
    to be the final full $2e$-state of $x$ w.r.t.\ $\{X_{i}\}_{i\leq
      2e}$ and $\{Y_{i}\}_{i\leq 2e}$, where $X_{2i} = S_{i}$,
    $X_{2i+1} =\breve{S}_{i}$, $Y_{2i}= T_{i}$, and $ Y_{2i+1} =
    \breve{T}_{i}$.
  \item Let $\nu^S(e,x)$ be the final $e$-split state of $x$ measured
    w.r.t.\ $\{W_i \searrow A\}_{i \leq e}$ and
    $\{\Gamma^{-1}(\tilde{S}_i)\}_{i \leq e}$.  Let
    $\widehat{\nu}^S(e,\xhat)$ be the final $e$-split state of $\xhat$
    measured w.r.t.\ $\{\Gamma(W_i \searrow A)\}_{i \leq e}$ and
    $\{\tilde{S}_i\}_{i \leq e}$.
  \item Let $\nu^S(\alpha,x)$ be the final $|\alpha|$-split state of
    $x$ measured w.r.t. $\{(W_i \searrow A)\}_{i \leq e}$ and
    $\{S_{\beta}\}_{\beta \subseteq \alpha}$.  (Careful---this is not
    the same as $\nu^S(|\alpha|,x)$.)
  \item Every $2e$-state is an \emph{$e$-split state} and $\nu =
    \langle 2e, \sigma, \tau\rangle$ is a \emph{reasonable $e$-split
      state} if for all $i \leq e$, exactly one of $2i$ or $2i+1$ is
    in $\sigma$, and exactly one of $2i$ or $2i+1$ is in $\tau$.
  \item For every $e$-split state $\nu$ and $\alpha$ such that
    $|\alpha| =e$, let
   \begin{multline*}
     D^A_{\nu,\alpha} = \{x : \exists s \text{ such that } x \in
     A_{s+1} - A_{s} \text{ and } \nu = \nu^S(e,x,s) \\ \text{ w.r.t.
     }  \{(W_i \searrow A)_s\}_{i \leq e,s<\omega}, \{(A \backslash
     W_i)_s\}_{i \leq e,s<\omega}, \\ \{S_{\beta,s}\}_{\beta \subseteq
       \alpha,s< \omega}, \text{ and } \{\breve{S}_{\beta,s}\}_{\beta
       \subseteq \alpha,s< \omega}.
  \end{multline*}
    \end{enumerate}
 
\end{definition}

Let $\nu$ be a reasonable $e$-split state.  Then $X_\nu= \{ x |
\nu^S(e,x) = \nu\}$ is a Boolean combination of splits of $A$ and
hence $X_\nu$ is also a split of $A$.  $\Gamma$ is an isomorphism
between $\S_\R(A)$ and $\S_\R(\Ahat)$ (modulo the computable subsets
of $A$).  Hence $\Gamma$ is an isomorphism between $\E_A \join
\Gamma^{-1}(\tilde{\B})$ and $\Gamma(\E_A) \join \tilde{\B}$ (again
modulo the computable subsets of $A$). Therefore, $X_\nu$ is
computable iff $\Gamma(X_\nu)$ is computable.  So, for all reasonable
$e$-split states $\nu$, $\{ x | \nu^S(e,x) = \nu\}$ is computable iff
$\{\xhat | \widehat{\nu}^S(e,\xhat)=\nu\}$ is computable.

Since $S_\alpha$ are entry sets, $x \in D^A_{\nu,\alpha}$ iff
$\nu^S(\alpha,x) = \nu$.  Therefore $\{x: \nu^S(\alpha,x) = \nu$ is
computable iff $D^A_{\nu,\alpha}$ is computable.

By Lemma~\ref{sec:preserve}, to show $\B$ is isomorphic via $\Theta$
to $\tilde{\B}$ it is enough to show, for all $\beta, \gamma$,
$S_\beta - S_\gamma$ is computable iff $\Theta(S_\beta) -
\Theta(S_\gamma) = \tilde{S}_{|\beta|} - \tilde{S}_{|\gamma|}$ is
computable.  Let $\alpha$ be the longer of $\beta$ and $\gamma$. Then
$$
S_\beta - S_\gamma = \bigsqcup \{D^A_{\nu,\alpha}: \nu = \langle
2|\alpha|,\sigma,\tau\rangle, 2|\beta|\in \tau, \text{ and } 2|\gamma|
\not\in \tau\}.$$
Therefore, it is more than enough to show, for all
reasonable $e$-split states $\nu$ and all $\alpha \subset f$, if
$|\alpha| = e$ then $D^A_{\nu,\alpha}$ is computable iff $\{ x |
\nu^S(e,x) = \nu\}$ is computable.

Hence from this point forward we will just work on constructing
$S_\alpha$ and $\breve{S}_\alpha$ such that for all reasonable
$e$-split states $\nu$ and all $\alpha \subset f$, if $|\alpha| = e$
then 
\begin{equation*}  \tag*{$\mathcal{R}_{\nu}$}
  D ^A_{\nu,\alpha} \text{ is computable iff }
   \{ x | \nu^S(e,x) = \nu\} \text{ is computable}.
\end{equation*}

(Let $\Theta(W_i \searrow A) = \Gamma(W_i \searrow A)$ and
$\Theta^{-1}(\Gamma(W_i \searrow A)) = W_i \searrow A$. Then almost the
same argument shows that $\Theta$ is an isomorphism between $\E_A$ and
$\Gamma(\E_A)$ and, in fact, $\E_A \join \B$ is isomorphic via $\Theta$ to
$\Gamma({\E}_A) \join \tilde{\B}$.)

If we succeed in meeting $\mathcal{R}_\alpha$ then $\Theta$ will be an
isomorphism as desired.  As we will see it turns out to do this it
enough to know for which for all reasonable $e$-splits states and
$\alpha$, $\{ x | \nu^S(e,x) = \nu\}$ is infinite.

Determining whether $\{ x | \nu^S(e,x) = \nu\}$ is infinite is
$\Delta^0_3$: Are there $i_k$ and $j_k$, for $k \leq e$, and
infinitely many $x$ and stages $s$ such that for all $k \leq e$, $
\Gamma^{-1}(\tilde{S}_k) = W_{i_k}$, $\breve{W}_{i_k} = W_{j_k}$, $x
\in W_{i_k,s} \sqcup W_{j_k,s}$, and $\nu^S(e,x,s) = \nu$, where
$\nu^S(e,x,s)$ is measured w.r.t.\ $\{(W_k \searrow A)_s\}_{k \leq
  e,s<\omega}$, $\{(A \backslash W_k)_s\}_{i \leq e,s<\omega}$,
$\{W_{i_k,s}\}_{k\leq e, s < \omega}$, and $\{W_{j_k,s}\}_{k\leq e, s
  < \omega}$.  Recall $\Gamma$ is $\Delta^0_3$ and since we know
$S=\Gamma^{-1}(\tilde{S}_k)$ is a split of $A$ we can find $\breve{S}$
using an oracle for $\mathbf{0''}$.  This also shows that $\{ x |
\nu^S(e,x) = \nu\}$ is a \ce set and a split of $A$.

Hence it is straightforward to construct a tree $Tr$, with a true path
$f$ and an approximation $f_s$ to $f$ such that $f = \liminf_s f_s$,
if $\alpha \in Tr$ then $\alpha$ is outfitted with a set of reasonable
$|\alpha|$-split states, $\mathcal{M}_{\alpha}$, and if $\alpha
\subset f$ then $\nu \in \mathcal{M}_{\alpha}$ iff $\{ x | \nu^S(e,x)
= \nu\}$ is infinite.  Furthermore we can assume that if $\beta
\subset \alpha$ and $\nu \in \mathcal{M}_{\alpha}$ then
$\nu\restriction 2|\beta| \in \mathcal{M}_{\beta}$ and that $|f_s|
=s$, for all $s$.  In the interest of space and energy we are not
going to go into the details.  Similar constructions with all the
details can be found in Section~\ref{tree2}, \citet{Cholak:95},
\citet{Cholak:94*1}, and \citet{Weber:04}.

Using the approximation to the true path we will construct a function
$\alpha(x,s)$ for all $x$ and $s$.  If $s < x$, let
$\alpha(x,s)\uparrow$.  Let $\alpha(x,x) = f_x$.  For $s \geq x$, if
$f_{s+1} <_L \alpha(x,s)$ then let $\alpha(x,s+1) = f_{s+1}$.

If $x$ enters $A$ at stage $s+1$ look for the greatest $\beta
\subseteq \alpha(x,s)$ where we can enumerate $x$ into $S_{\gamma,s}$
and $\breve{S}_{\gamma,s}$, for $\gamma \subseteq \beta$, such that
$\nu^S(\beta,x,s) = \nu \in \mathcal{M}_{\beta}$ and, for all $\beta'
\subset \beta$, if we can enumerate $x$ into $S_{\gamma,s}$ and
$\breve{S}_{\gamma,s}$, for $\gamma \subseteq \beta'$, such that
$\nu^S(\beta,x,s) = \nu' \in \mathcal{M}_{\beta'}$ then
$D^A_{\nu',\beta',s} \neq \emptyset$. If there are several possible
$\nu$, arbitrarily choose the one where $D^A_{\nu,\beta,s}$ is the
smallest.  Enumerate $x$ such that $\nu^S(\beta,x,s) = \nu$.  For all
$\gamma$, if $\gamma \nsubseteq \beta$ or $\beta$ does not exist, add
$x$ to $\breve{S}_{\gamma,s}$.

For any $\beta \subset f$, let $s_{\beta}$ be such that if $f_t <_L
\beta$ then $t< s_{\beta}$, if $\{ x | \nu^S(|\beta|,x) = \nu'\}$ is
finite and $\nu^S(|\beta|,x) = \nu'$ then $x<s_{\beta}$, and if $\nu
\in \mathcal{M}_{\beta}$ then $D^A_{\nu,\beta,s} \neq \emptyset$ (by
induction on $\beta$ it is not hard to show that such a stage exists).  For
each $x\geq s_{\beta}$ we can effectively find a stage $s_{\beta,x}$
such that for all $s' \geq s_{\beta,x}$, $\beta \subseteq
\alpha(x,s')$.  Let $R_\beta$ be the set of $x$ such that either $x <
s_{\beta}$ and $x \in A$ or $x\geq s_{\beta}$ and $x \in
A_{s_{\beta,x}}$.  $R_\beta$ is a computable subset of $A$.

\begin{lemma}
  If $\alpha \subset f$ and $\nu$ is a reasonable $|\alpha|$-split
  state then $D^A_{\nu,\alpha}$ is computable iff $\{ x |
  \nu^S(|\alpha|,x) = \nu\}$ is computable.
\end{lemma}

\begin{proof}
  Let $|\alpha| =e$ and $\nu = \langle 2e, \sigma, \tau\rangle$.
  
  ($\Rightarrow$) Assume $\{ x | \nu^S(e,x) = \nu\}$ is not
  computable.  We must show $D^A_{\nu,\alpha}$ is not computable.
  Assume otherwise.  Hence there is an $i > e$ such that $W_i =
  D^A_{\nu,\alpha}$, and $A \searrow W_i =\emptyset$.  There must
  exist a reasonable $i$-split state $\nu' = \langle 2i , \sigma',
  \tau' \rangle$ such that $\sigma' \restriction 2e = \sigma$, $2i+1
  \in \sigma'$, $\tau' \restriction 2e = \tau$, and $\{ x | \nu^S(i,x)
  = \nu'\}$ is not computable. (Otherwise $\{ x | \nu^S(e,x) = \nu\}$
  is computably contained in a computable set, $W_i$, and hence is
  computable.)  Therefore $\{ x | \nu^S(i,x) = \nu'\} - R_\beta$ is
  infinite. Hence, by the above construction, there is an $x$ such
  that $x \in D^A_{\nu',\beta}$.  This same $x$ is in
  $D^A_{\nu,\alpha}$ but not in $W_i$.  Contradiction.
 
  ($\Leftarrow$) Assume $\{ x | \nu^S(e,x) = \nu\}$ is computable.
  Hence there is an $i > e$ such that $W_i = \{ x | \nu^S(e,x) =
  \nu\}$, and $A \searrow W_i =\emptyset$.  Let $\beta \subset f$ and
  $|\beta| =i$.  For $j \geq i$, if $\nu' = \langle 2j , \sigma',
  \tau' \rangle$, $\sigma' \restriction 2e = \sigma$, $2i+1 \in
  \sigma'$, and $\tau' \restriction 2e = \tau$, then $\{ x | \nu^S(j,x)
  = \nu'\}$ is not infinite.  Hence for all $\gamma \supseteq \beta$,
  $\nu' \not\in \mathcal{M}_{\gamma}$.  Let $x \in \overline{W_i} -
  R_\beta$ enter $A$ at stage $s+1$.  Then $\nu^S(i,x) = \nu' \in
  \mathcal{M}_{\beta}$ and $\nu' \restriction 2e \neq \nu$.  Hence, by
  the above construction $\nu^S(\beta,x,s) \neq \nu$.  Therefore if
  $x$ enters $A$ at stage $s+1$ and $\nu^S(\beta,x,s) = \nu$ then $x
  \in R_\beta$ or $x \in W_i$.  Thus $D^A_{\nu,\alpha}$ is computable.
\end{proof}

Therefore $\Theta$ is an isomorphism between $\B$ and $\widehat{\B}$.
Thus (1) holds.  The next lemma proves (2).

\begin{lemma}
  If $\alpha \subset f$, $|\alpha| = e$, and $\tilde{S}_e$ supports
  $\What$, then $S_\alpha$ supports $X= \Psi^{-1}(\widehat{W})$.
\end{lemma}

\begin{proof}
  Since $\Psi$ is an automorphism of $\E^*$ taking $A$ to $\Ahat$,
  $\Psi^{-1}(\tilde{S}_e)$ supports $X$.  Since
  $\Gamma$ is induced by $\Psi$, $\Gamma^{-1}(\tilde{S}_e)$ supports
  $X$.  Let $i > e$ such that $W_i = (X -A) \sqcup
  \Gamma^{-1}(\tilde{S}_e)$.  Hence $W_i \searrow A$ supports $X$.  If
  $(W_i \searrow A) \subseteq_\R Y$ then $Y$ supports $X$.  Hence it
  is enough to show $(W_i \searrow A) \subseteq_\R S_\alpha$.
  
  Let $\beta \subset f$ such that $|\beta| =i$.  For $j \geq i$, if
  $\nu = \langle 2j , \sigma, \tau \rangle$, $2i \in \sigma$, and $\{
  x | \nu^S(j,x) = \nu\}$ is infinite then $2e \in \tau$.  Hence for
  all $\gamma \supseteq \beta$, if $\nu = \langle |\gamma|, \sigma,
  \tau\rangle \in \mathcal{M}_{\gamma}$ and $2i \in \sigma$ then $2e
  \in \tau$.  Let $x \in {W_i} - R_\beta$ enter $A$ at stage $s+1$.
  Then $\nu^S(i,x) = \nu \in \mathcal{M}_{\beta}$.  Hence, by the
  above construction, for almost all such $x$, $x \in S_\alpha$.
  Hence $(W_i \searrow A) \subseteq^* S_\alpha \cup R_\beta$. 
\end{proof}

\section{A definable orbit which is not a $\Delta^0_3$
  orbit} \label{sorbits}
  
For $\E^*$, all the previously known orbits are actually orbits under
$\Delta^0_3$-automorphisms.  And a good number of those are also
definable in the sense that there is an elementary formula, $\varphi(X)$,
in the language of $\E^*$ such that $\varphi(A)$ iff $A$ is in the
orbit under question.  Examples include maximal sets, creative sets,
hemimaximal sets, and quasi-maximal sets.

The following is a definable orbit $\mathcal{O}$, which is not a
$\Delta^0_3$ orbit.  It is the first example of an orbit which is not
an orbit under $\Delta^0_3$-automorphisms.  It is an orbit under
$\Delta^0_5$-automorphisms.

In the mid 1990s, Cholak and Downey incorrectly claimed to construct
a pair of $\Delta^0_4$-automorphic \ce sets which were not
$\Delta^0_3$-automorphic. In addition, we show this claim is correct
by showing there are two such sets in $\mathcal{O}$.

\subsection{The orbit $\mathcal{O}$}

Assume that $A$ is not computable.

\begin{definition}
  $F$ is \emph{$A$-special} if $F$ is not computable, $F \cap A =
  \emptyset$, and, for all $V$, if $V \cap A = \emptyset$ then $V- F$
  is computably enumerable.
\end{definition}

\begin{lemma}\label{sec:orbit-mathcalo}
  Assume $F_0$ and $F_1$ are $A$-special sets and $R$ is computable
  set disjoint from $A$.
  \begin{enumerate}
  \item Either $F_1 - F_0$ is computable or $A$-special.
  \item  If $F_0 \cap R = \emptyset$ then $F_0 \sqcup R$ is $A$-special.
  \item If $F_0 \cap F_1 = \emptyset$ then $F_0 \sqcup F_1$ is
    $A$-special.
  \item $F_0 \cup F_1$ is $A$-special.
  \item If $W \subseteq R$ then $W$ is not $A$-special.
  
  \end{enumerate}
\end{lemma}

\begin{proof}
  (1) $V - (F_1 - F_0) = (V -F_1) \cup (V \cap F_0)$. So if $F_1-F_0$
  is not computable, it is $A$-special.
  
  (2) $V - (F_0 \sqcup R) = (V-F_0)-R$. If $F_0 \sqcup R$ is
  computable then $F_0$ is computable.
  
  (3) $ V - (F_0 \cup F_1) = (V - F_0) -F_1$. If $F_0 \sqcup F_1$ is
  computable then $F_0$ is computable.
  
  (4) $F_0 \cup F_1 = F_0 \sqcup (F_1 - F_0)$.  Now apply (1) followed
  by one of (2) or (3).
  
  (5) If for all $V$, if $V \cap A = \emptyset$ then $V- W$ is
  computably enumerable then $\overline{W} = (R-W) \cup \overline{R}$.
\end{proof}

\begin{definition}\label{sec:orbit-mathcalo-1}
  Let $\varphi(A)$ be the conjunction of the following 3 statements:
  \begin{enumerate}
  \item $\forall F$ if $F$ is $A$-special then $\exists G$ such that
    $G$ is $A$-special and $F \cap G = \emptyset$;
  \item $\forall W$ if $W \cap A = \emptyset$ then $\exists F$ such
    that $F$ is $A$-special and $W \subseteq^* F$;
  \item $\forall W \exists F$ such that $F \cap A = \emptyset$ and
    either $W \subseteq^* F \sqcup A$ or $W \cup F \cup A =^* \omega$.
  \end{enumerate}
\end{definition}

\begin{definition}
  A list of \ce sets, $\mathcal{F} = \{ F_i : i \in \omega\}$, is an
  \emph{$A$-special list} iff $\mathcal{F}$ is a list of pairwise
  disjoint noncomputable sets, $F_0 = A$, and for all $W$ there is an
  $i$ such that $W \subseteq^* \bigsqcup_{l\leq i}F_l$ or $W \cup
  \bigsqcup_{l\leq i}F_l =^* \omega$.  We say that $\mathcal{F}$ is a
  \emph{$\Gamma$ $A$-special list} if $\mathcal{F}$ is an $A$-special
  list and there is a function $f$ with property $\Gamma$ such that
  $F_i = W_{f(i)}$.
\end{definition}

Note that for any $i$, $\bigsqcup_{ l \geq i}F_l$ is not \ce and hence
there cannot be an effective $A$-special list.  The automorphic image
under $\Phi$ of an $A$-special list is a $\Phi(A)$-special list.

\begin{lemma}\label{sec:orbit-mathcalo-5}
  Assume that an $A$-special list exists and that $V \cap A =
  \emptyset$.  Then $V \subseteq^* \bigsqcup_{0< l \leq i} F_l$, for
  some $i$.
\end{lemma}

\begin{proof}
  If $V \cup \bigsqcup_{l\leq i}F_l =^* \omega$, for some $i$, then
  $(V \cup \bigsqcup_{0< l\leq i}F_l) \sqcup A =^* \omega$ and hence
  $A$ is computable. Contradiction.
\end{proof}

\begin{lemma}\label{sec:orbit-mathcalo-3}
  $\varphi(A)$ iff an $\mathbf{0^{(4)}}$ $A$-special list exists.
\end{lemma}

\begin{proof}
  ($\Rightarrow$) Let $F_0 = A$.  Assume, by induction, for $0< j <
  i$, that $F_j$ are defined such that they are pairwise disjoint,
  $A$-special, either $W_j \subseteq^* \bigsqcup_{l\leq j}F_l$ or $W_j
  \cup \bigsqcup_{l\leq j}F_l =^* \omega$, and $\bigsqcup_{0< j < i}
  F_j$ is $A$-special.  Since $\phi(A)$ holds, the third clause of
  Definition~\ref{sec:orbit-mathcalo-1} holds for $W_i$ and hence
  there is an $F$ such that $F \cap A = \emptyset$ and either $W_i
  \subseteq^* F \sqcup A$ or $W_i \cup F \cup A =^* \omega$.  By the
  second clause of Definition~\ref{sec:orbit-mathcalo-1} and the fact
  that $A$-special sets are disjoint from $A$, we can assume $F$ is
  $A$-special.  Hence, by Lemma~\ref{sec:orbit-mathcalo},
  $\bigsqcup_{j < i} F_j \cup F$ is $A$-special and $F - \bigsqcup_{j
    < i} F_j$ is either computable or $A$-special.  If $F -
  \bigsqcup_{j < i} F_j$ is $A$-special let $F_i = F - \bigsqcup_{j <
    i} F_j$.  Otherwise apply the first clause of
  Definition~\ref{sec:orbit-mathcalo-1} to $\bigsqcup_{j < i} F_j \cup
  F$ to get an $A$-special $G$ and let $F_i = G \sqcup (F -
  \bigsqcup_{j < i} F_j)$ which is $A$-special by
  Lemma~\ref{sec:orbit-mathcalo}.  Again by
  Lemma~\ref{sec:orbit-mathcalo}, $\bigsqcup_{0< j \leq i} F_j$ is
  $A$-special.
  
  If $X$ and $Y$ are \ce sets then whether $Y - X$ is \ce is
  $\Sigma^0_3$.  So whether $F$ is $A$-special is $\Pi^0_4$.  Since
  $\varphi(A)$ holds, given $W$, there exists an $A$-special set $F$
  such that either $W \subseteq^* F \cup A$ or $W \cup F \cup A =^*
  \omega$.  Hence we can try all possible $F$ using $\bf{0^{(4)}}$ to
  test if the $F$ being considered has the correct properties. Since
  such an $F$ exists this algorithm will converge and is computable in
  $\bf{0^{(4)}}$.  Going from $F$ to $F_i$ is also computable in
  $\bf{0^{(4)}}$. Hence the $A$-special list constructed is computable
  in $\bf{0^{(4)}}$.
  
  ($\Leftarrow$) By Lemma~\ref{sec:orbit-mathcalo}, it is enough to
  show that for all $j \geq 1$, $F_j$ is $A$-special.  To show $F_j$
  is $A$-special it is enough to show that if $V \cap A =\emptyset$ then
  $V-F_j$ is computably enumerable. Assume $V \cap A =\emptyset$.
  Then, by Lemma~\ref{sec:orbit-mathcalo-5}, $V \subseteq^*
  \bigsqcup_{0< l \leq i} F_l$, for some $i$.  So $V - F_j =^* V \cap
  \bigsqcup_{0< l \leq i \wedge l \neq j} F_l$ is a \ce set.
\end{proof}

\begin{theorem}\label{sec:orbit-mathcalo-2}
  Given an $\bf{a}$ $A$-special list, $\mathcal{F}$, and an
  $\widehat{\bf{a}}$ $\Ahat$-special list, $\widehat{\mathcal{F}}$,
  there is a $\bf{0''} \join \bf{a} \join$
  $\widehat{\bf{a}}$-automorphism $\Theta$ of $\E^*$ taking $A$ to
  $\Ahat$.
\end{theorem}

\begin{proof}
  By Theorem~\ref{apcomputable}, there is an isomorphism $\Theta_i$
  between $F_i$ to $\widehat{F}_i$ preserving computable sets.  Given
  $W_e$ define $\Theta(W_e)$ as follows: If $W_e \subseteq^*
  \bigsqcup_{l\leq i}F_l$ then $\Theta(W_e) = \bigsqcup_{l \leq
    i}\Theta_l(W_e \cap F_l)$.  Otherwise there is a computable set
  $R$ such that $R \subseteq^* \bigsqcup_{l \leq i} F_l$ and $R \cup W_e
  =^* \omega$.  For all $l \leq i$, $R \cap F_l$ is computable.
  Therefore, since $\Theta_l$ preserves computable sets, $\Theta(R) =
  \bigsqcup_{l \leq i}\Theta_l(R \cap F_l)$ is computable.  Let
  $$\Theta(W_e) = \overline{\Theta(R)} \sqcup \bigsqcup_{l \leq
    i}\Theta_l(W_e \cap R \cap F_l).$$
  
  $\Theta$ is an automorphism of $\E^*$ such that $\Theta(A) = \Ahat$.
  By Theorem~\ref{apcomputable}, an index for $\Theta_i$ can be found
  uniformly from indices for $F_i$ and $\widehat{F}_i$.  The remaining
  division into cases can be done using a $\bf{0''}$ oracle.
\end{proof}

\begin{theorem}\label{orbit1}
  The collection of $A$ such that $\phi(A)$ forms a $\Delta^0_5$
  orbit $\mathcal{O}$.
\end{theorem}
  
\begin{proof}
  This follows from Theorems~\ref{sec:orbit-mathcalo-3} and
  \ref{sec:orbit-mathcalo-2}.
\end{proof}

\begin{corollary}
  If $\mathcal{F}$ is an $A$-special list then, for all $i$, $F_i$ is
  automorphic to $A$.
\end{corollary}

\begin{proof}
  The list formed by switching $F_i$ and $A$ is an $F_i$-special list. 
\end{proof}
 
\subsubsection{$\mathcal{O}$ is not a $\Delta^0_3$ orbit}
\label{sec:orbit-mathcalo-4}

\begin{theorem} \label{orbit2}
  There are \ce sets $A$ and $\Ahat$ such that $\phi(A)$ and
  $\phi(\Ahat)$, $A$ and $\Ahat$ are $\Delta^0_4$-automorphic but not
  $\Delta^0_3$-automorphic.
\end{theorem}

This theorem follows from the next two lemmas.  

\begin{lemma} \label{Aexists} 
  There exists $A$ such that a $\bf{0''}$ $A$-special list
  $\mathcal{F}$ exists.
\end{lemma}

\begin{lemma} \label{Ahatexists}
  There exists $\Ahat$ such that a $\bf{0'''}$ $\Ahat$-special list
  $\widehat{\mathcal{F}}$ exists but no $\bf{0''}$ $\Ahat$-special
  list exists.
\end{lemma}

The proofs of these lemmas follow in
Section~\ref{sec:proofs-lemm-refa}.

\begin{proof}[Proof of Theorem~\ref{orbit2} 
  from Lemmas~\ref{Aexists} and \ref{Ahatexists}] Assume that
  $\mathcal{F}$ is the $A$-special list given by Lemma~\ref{Aexists}
  and $\widehat{\mathcal{F}}$ is the $\Ahat$-special list given by
  Lemma~\ref{Ahatexists}.  By Lemma~\ref{sec:orbit-mathcalo-2}, $A$
  and $\Ahat$ are in $\mathcal{O}$ and are $\Delta^0_4$-automorphic.
  
  Let $f$ witness that $\mathcal{F}$ is a $\bf{0''}$ $A$-special list.
  Assume that $A$ and $\Ahat$ are $\Delta^0_3$ automorphic via
  $\Phi(W_e) = W_{g(e)}$ then $\{ W_{g((f(i))} | i \in \omega\}$ is
  $\bf{0''}$ $\Ahat$-special list.  Therefore $A$ and $\Ahat$ cannot
  be in the same $\Delta^0_3$ orbit.
\end{proof}

The following lemma and corollary are needed for the proof of
Lemma~\ref{Ahatexists}.

\begin{lemma}\label{sec:mathc-not-delt}
  If a $\bf{0''}$ $A$-special list $\mathcal{F} = \{ F_i : i \in
  \omega\}$ exists then there is a function $d$ computable in
  $\bf{0''}$ such that if $W_e \cap A = \emptyset$ then $W_{d(e)} \cap
  (W_e \cup A) =\emptyset$ and $W_{d(e)}$ is $A$-special.
\end{lemma}

\begin{proof}
  If $W_e \cap A \neq \emptyset$ (whether this occurs is computable in
  $\bf{0''}$) then let $d(e) = 0$.  Assume $W_e \cap A =\emptyset$.
  Let $f$ witness that $\mathcal{F}$ is $\bf{0''}$.  Then, by
  Lemma~\ref{sec:orbit-mathcalo-5}, $W_e \subseteq^* \bigsqcup_{0< l
    \leq i} F_l$, for some $i$. Using $f$, the least such $i$ can be
  found computably in $\bf{0''}$.  Let $d(e) = f(i+1)$.
\end{proof}

\begin{corollary}\label{sec:mathc-not-delt-1}
  Assume for all $e$, there are $e'$ and $d$ such that $W_{e'} \cap A
  = \emptyset$ and if $W_{\phi(\langle e', d \rangle)}$ is cofinite
  then either $W_d \cap (W_{e'} \cup A) \neq \emptyset$ or $W_d$ is
  not $A$-special.  Then $A$ does not have a $\bf{0''}$ $A $-special
  list.
\end{corollary}

\begin{proof}
  Assume $A$ has a $\bf{0''}$ $A$-special list.  Apply
  Lemma~\ref{sec:mathc-not-delt} to get ${g}$.  The graph of $g$ is a
  $\Delta^0_3$ set and hence a $\Sigma^0_3$ set. $\Cof$ is
  $\Sigma^0_3$-complete.  Hence there is an $e$ such that, for all
  $e'$, $W_{\phi(\langle e', d(e') \rangle)}$ is cofinite and if
  $W_{e'} \cap A =\emptyset$ then $W_{d(e')} \cap (W_{e'} \cup A) =
  \emptyset$ and $W_{d(e')}$ is $A$-special.  Furthermore, since we
  are reducing the graph of a function to $\Cof$, for all $e'$, if $d
  \neq d(e')$ then $W_{\phi(\langle e', d \rangle)}$ is not cofinite.
  Contradiction.
\end{proof}

\subsection{Proofs of Lemmas~\ref{Aexists} and \ref{Ahatexists}}
\label{sec:proofs-lemm-refa}

First we will focus on Lemma~\ref{Aexists}.  Rather than focusing on
$A$ we will first focus on constructing the $A$-special list
$\mathcal{F}$.  This will be a tree argument and very similar to the
$\Delta^0_3$-isomorphism method.  At each node $\alpha \in T$ we will
build a \ce set, $F_{\alpha}$.  The goal is to build the
$F_{\alpha}$s such that if $F_i = F_{|\alpha|}$, for $\alpha \subset
f$, where $f$ is the true path, then $\mathcal{F} = \{ F_i : i \in
\omega\}$ is an $A$-special list.

\subsubsection{The requirements}

We will construct the $F_{\alpha}$s as pairwise disjoint 
noncomputable
sets, for $\alpha \subset f$. $F_{\alpha}$ must be noncomputable.
Hence we must meet the following requirements for all $\alpha \subset
f$ and all $e$:
\begin{equation*}
  \tag*{$\mathcal{R}_{\alpha,e}$:}
     \overline{F}_{\alpha} \neq W_e.
\end{equation*}
In addition, we will meet the following requirement for all $\alpha
\subset f$:
\begin{equation*}
  \tag*{$\mathcal{N}_{\alpha}$:}
     \text{either } W_{|\alpha|} \subseteq^*  
   \bigsqcup_{\beta \subseteq \alpha} F_{\beta} \text{ or }
     W_{|\alpha|} \cup \bigsqcup_{\beta \subseteq \alpha} 
    F_{\beta} =^* \omega.
\end{equation*}
Before we can discuss how we will meet these requirements we need the
following remark.

\begin{remark}[The position function $\alpha(x,s)$] \label{standard}
  Given the approximation to the true path at stage $s$, $f_s$, we
  will determine the position function $\alpha(x,s)$ by the following
  rules: $x$ is \emph{$\alpha$-legal at stage $s$} if $\alpha(x,s-1) =
  \alpha^-$ (recall $\alpha^-$ is the node before $\alpha$ in the
  tree), $x$ is $\alpha^-$-allowed (defined below) and for all stages
  $t$, if $x \leq t \leq s$, then $\alpha \leq_L f_t$.  If $\alpha
  \subseteq f_s$ and $x$ is $\alpha$-legal then let $\alpha(x,s) =
  \alpha$ (move $x$ downward into $\alpha$). If $f_s <_L
  \alpha(x,s-1)$ then let $\alpha(x,s) = \alpha(x,s-1) \cap f_s$.
\end{remark}

\subsubsection{Action for $\mathcal{R}_{\alpha,e}$}

Meeting $\mathcal{R}_{\alpha,e}$ is straightforward.  But we are going
to break it into parts, ensuring that there are possible witnesses and
actually taking action to meet $\mathcal{R}_{\alpha,e}$.  

\emph{Getting witnesses:} For each $\beta$ and each stage $s$, we will
pick a $x_{\beta,s}$.  We will hold $x_{\beta,s}$ out of all
$F_{\gamma}$, for $\gamma \supset\beta$ but allow $x_{\beta,s}$ to
possibly enter $F_{\gamma}$, for $\gamma \subseteq \beta$.  If
$x_{\beta,s}$ enters some $F_{\gamma}$ at stage $s$ (or does not
exist yet), then, at the next stage $t$, such that $\beta \subseteq
f_t$ and there is an $x$ with $\alpha(x,s) = \beta$ and $x \notin
\bigsqcup_{\beta \in T} F_{\beta,s}$, we will choose the least such
$x$ as $x_{\beta,t}$; until that stage $t$, $x_{\beta,s}$ does not
exist. Otherwise $x_{\beta,s}$ remains the same from stage to stage.

\emph{Placing witnesses into $F_\alpha$:} Now if $\alpha \subseteq
f_s$, $W_{e,s} \cap F_{\alpha,s} = \emptyset$, $|\alpha| \leq e$, and
there is an $x$ where $|\alpha(x,s)| \geq |\alpha| + e$, $x \in
W_{e,s}$ and $x \notin \bigsqcup_{\beta \in T} F_{\beta}$, then add
$x$ to $F_{\alpha}$ at stage $s$.

Assume that for all $\gamma \subset f$, $x_{\gamma} = \lim_s
x_{\gamma,s}$ exists.  Then if $\overline{F}_{\alpha} = W_e$ then it
is straightforward to show that at some stage $s$ we will add an $x$ to
$F_{\alpha}$ to meet $\mathcal{R}_{\alpha,e}$.

Notice that only finitely many $\mathcal{R}_{\alpha,e}$ are possibly
interested in $x_{\gamma,s}$.  So if $x_{\gamma}$ fails to exist it
is not due to action for $\mathcal{R}_{\alpha,e}$ but action for some
$\mathcal{N}_{\beta}$.

\subsubsection{Action for $\mathcal{N}_{\alpha}$}

We will meet $\mathcal{N}_{\alpha}$ as follows: First of all no action
is taken at stage $s$ if $x_{\alpha,s}$ does not exist.  Furthermore,
we never $\alpha$-allow $x_{\alpha,s}$. Otherwise the desired action
at $\alpha$ breaks into cases depending on whether $W_{\alpha}$ is
infinite or not, where
\begin{equation*}
  W_{\alpha} = \{ x | \exists s (\alpha^- \subseteq \alpha(x,s) \wedge
   x \text{ is $\alpha^-$-allowed }
  \wedge x \in W_{|\alpha|,s})\}.  %
\end{equation*}
If $\alpha$ believes $W_{\alpha}$ is finite we \emph{$\alpha$-allow}
half of the balls which arrive at $\alpha$ (hence these balls can move
downward) and put all but one ball, $x_{\alpha,s}$, of the other half
into $F_{\alpha}$ (like $x_{\alpha,s}$, the balls in $F_{\alpha}$, are
never $\alpha$-allowed).  Assume $\alpha$ believes $W_{\alpha}$ is
infinite.  Half of the balls which arrive at $\alpha$ in $W_{\alpha}$
will be $\alpha$-allowed immediately.  Otherwise if $\alpha(x,s) =
\alpha$ and there have been $x$ many balls $\alpha$-allowed, we will
place $x$ into $F_{\alpha}$.

\subsubsection{The Verification}

Assume that for all $\alpha \subset f$, infinitely many balls are
$\alpha$-allowed (we will show this later).  Then, by induction on
$\alpha \subset f$, it is straightforward to show that $x_{\alpha}$
exists and hence $\mathcal{R}_{\beta,e}$ is met for $\beta \subset f$
and all $e$. And, again by induction on $\alpha \subset f$, is
straightforward to show, using the standard facts about $f_s$ and
$\alpha(x,s)$ and the above assumption, for almost all $x \notin
\bigsqcup_{\beta\subset\alpha} F_{\beta}$, there is a stage such that
either $x$ enters $F_{\alpha}$ or $x$ is $\alpha$-allowed.  Hence if
$W_{\alpha}$ is finite then $W_{|\alpha|} \subseteq^*
\bigsqcup_{\beta\subset\alpha} F_{\beta}$ and otherwise $W_{|\alpha|}
\cup \bigsqcup_{\beta\subseteq\alpha} F_{\beta} =^* \omega$.
Therefore, under the above assumption, $\mathcal{N}_{\alpha}$ is met.

Now we will show, by induction on $\alpha \subset f$, that infinitely
many balls are $\alpha$-allowed.  Assume this is true for $\alpha^-$.
Almost all of the balls which are $\alpha^-$-allowed will arrive at
$\alpha$ at some later stage (i.e., there is a stage $t$ such that
$\alpha \subseteq \alpha(x,t)$).  Hence at almost all stages,
$x_{\alpha,s}$ exists.  Therefore if $W_{\alpha}$ is finite then half
of those balls which arrive at $\alpha$ will be $\alpha$-allowed.  If
$W_{\alpha}$ is infinite then infinitely many balls arrive at $\alpha$
in $W_{\alpha}$, half of which are $\alpha$-allowed.

Hence the only thing needed to complete the proof of
Lemma~\ref{Aexists} is to construct the tree $T$, the true path $f$,
and the approximation to the true path at stage $s$, $f_s$.  But since
we want to use the same tree and related materials for the proof of
Lemma~\ref{Ahatexists}, we will delay this until Section~\ref{tree2}.

\subsubsection{Changes needed for the proof of Lemma~\ref{Ahatexists}}
\label{sec:changes-needed-proof}  

Rather than proving Lemma~\ref{Ahatexists} we will prove its unhatted
dual. We are going to make use of Lemma~\ref{sec:mathc-not-delt-1}.
We must meet the requirements:

\begin{equation*}
  \tag*{$\mathcal{Q}_{e}$:}
   \begin{split}  
     \text{there are }e' \text{ and } d \text{ such that }W_{e'} \cap
     A = \emptyset \text{ and if $W_{\phi(\langle e', d \rangle)}$ is
       cofinite} \\ \text{then either $W_d \cap (W_{e'} \cup A) \neq
       \emptyset$ or $W_d$ is not $A$-special.}\end{split}
\end{equation*}

By the Recursion Theorem we can assume there are computable functions
$g$ and $h$ such that $W_{g(\alpha)} = F_\alpha$ and $W_{h(\alpha)} =
\bigcup_{\lambda \subset \beta \subseteq \alpha} F_{\beta}$, for all
$\alpha \in T$ and $\alpha \neq \lambda$.  Recall $\lambda$ is the
empty node and $F_{\lambda} = A$. For all $\alpha \neq \lambda$,
$W_{h(\alpha)} \cap A = \emptyset$.

Assume that $\alpha$ is assigned to meet $\mathcal{Q}_e$.  $\alpha$
will use $W_{h(\alpha)}$ as $W_{e'}$.  We want to look for the least
$d$ and $l$ such that $[l, \infty) \subseteq W_{\phi_e(\langle
  h(\alpha),d\rangle)}$. We will use the tree to find $k$ and $l$ and
to assign $\alpha$ to $\mathcal{Q}_e$.

We will define the tree such that there are $d, l$ where $[l, \infty)
\subseteq W_{\phi_e(\langle h(\alpha),d\rangle)}$ iff there is a
unique $\beta$ such that $\alpha \subset \beta \subset f$ and $\beta$
believes there are $d, l < |\beta|$ such that $[l, \infty) \subseteq
W_{\phi_e(\langle h(\alpha),d\rangle)}$.  We will assume that the
$\mathcal{Q}_i$ are assigned in increasing order modulo finite injury
along the true path.  The finite injury along the true path will be
discussed below.

Assume that $\beta$ believes there are $d, l < |\beta|$ such that $[l,
\infty) \subseteq W_{\phi_e(\langle h(\alpha),d\rangle)}$. Since $d <
|\beta|$ there is a $\gamma \subset \beta$ with $|\gamma| =d$.
Furthermore, since we will continue to meet $\mathcal{N}_{\gamma}$,
either $W_{d} \subseteq^* \bigsqcup_{\delta \subseteq \gamma}
F_{\delta}$ or $W_{d} \cup \bigsqcup_{\delta \subseteq \gamma}
F_{\delta} =^* \omega$.  By Lemma~\ref{sec:orbit-mathcalo-5}, if
$W_{d} \cup \bigsqcup_{\delta \subseteq \gamma} F_{\delta} =^* \omega$
then $W_d \cap A \neq \emptyset$ and we have met $\mathcal{Q}_e$.  If
$W_{d} \subseteq^* \bigsqcup_{\delta \subseteq \alpha} F_{\delta}$
then we have met $\mathcal{Q}_e$.  Hence the only case where we must
take action to meet $\mathcal{Q}_e$ is when $W_{d} \subseteq^*
\bigsqcup_{\alpha \subset \delta \subseteq \gamma} F_{\delta}$. In
this case we will force $\bigsqcup_{\alpha \subset \delta \subseteq
  \gamma} F_{\delta}$ to be computable and hence, by
Lemma~\ref{sec:orbit-mathcalo} (5), $W_d$ is not $A$-special.  This
means we will have to later reconsider how we form the $A$-special
list.

Assume that $\beta$ must take action to meet $\mathcal{Q}_e$. $\beta$
will take action by changing how we meet $\mathcal{R}_{\gamma,e}$, for
all $\alpha \subseteq \gamma \subseteq \beta$.  Let $\alpha \subseteq
\gamma \subseteq \beta$.  The action taken for
$\mathcal{R}_{\gamma,e}$ is revised as follows: if $\gamma \subseteq
f_s$, $W_{e,s} \cap F_{\gamma,s} = \emptyset$, and there is an $x$ such
that $\beta \nsubseteq \alpha(x,s)$, $|\alpha(x,s)| \geq |\gamma| +
e$, $x \in W_{e,s}$ and $x \notin \bigsqcup_{\delta\in T} F_{\delta}$,
then add $x$ to $F_{\gamma}$ at stage $s$.  Now to help with the
creation of an $A$-special list we must injure all $\mathcal{Q}_i$
assigned to some $\gamma$ between $\alpha$ and $\beta$.  We will
assign then in increasing order to some $\delta$ where $\beta \subset
\delta$. This is finite injury along the true path.

If no $\alpha \subset \beta \subset f$ believes that it must take
action to meet $\mathcal{Q}_e$ then the above argument for the
verification of $\mathcal{R}_{\gamma,e}$ still holds and $F_{\gamma}$
is not computable.

Assume that some $\beta \subset f$ believes that it must take action
to meet $\mathcal{Q}_e$.  From the above verification, we know that
almost all $x$ either enter $\bigsqcup_{\delta \subseteq \beta}
F_{\delta}$ or are $\beta$-allowed.  By the above modification of the
action for $\mathcal{R}_{\gamma,e}$ once a ball either enters
$\bigsqcup_{\delta \subseteq \beta} F_{\delta}$ or is $\beta$-allowed
it cannot be used to meet $\mathcal{R}_{\gamma,e}$.  Hence $F_{\gamma}$
is computable and $\mathcal{Q}_e$ is met. 

The issue of an $A$-special list remains. Using the true path $f$ and
$\bf{0'''}$ we will inductively show how to construct an $A$-special
list.  Assume that we have built the list up to $i$ and have used
$\alpha_i \subset f$.  Let $\alpha^+$ be such that $\alpha \subset
\alpha^+ \subset f$ and $|\alpha^+| = |\alpha|+1$. Assume that
$\mathcal{Q}_e$ is assigned to $\alpha^+$ and by induction
$\mathcal{Q}_e$ is not injured from below.  Use $\bf{0'''}$ to see if
some $\beta \subset f$ takes action to meet $\mathcal{Q}_e$. If no
$\beta \subset f$ must take action to meet $\mathcal{Q}_e$ then
$F_{i+1} = F_{\alpha^+}$ is not computable and let $\alpha_{i+1} =
\alpha^+$.  Otherwise there is a $\beta \subset f$ which takes action
to meet $\mathcal{Q}_e$. In this case $F_\beta$ is not computable and
let $F_{i+1} = \bigsqcup_{\alpha_i \subset \gamma \subseteq \beta}
F_\gamma$ and $\alpha_{i+1} = \beta$.  In either case there is no
injury from below above $\alpha_{i+1}$.

\subsubsection{The tree $T$ and related definitions}\label{tree2}

We will define one tree which can be used for both lemmas.  We will
define $T$, the true path $f$, and the approximation to the true path
at stage $s$, $f_s$ via induction on the length of $\gamma$. 

We have to code a few items into $T$. At a node $\beta$ we must code
whether $W_\beta$ is infinite and whether there exists an $\alpha
\subset \beta$ and $e,d,l,s < \beta$ such that $\mathcal{Q}_e$ is
assigned by $\alpha$, $\alpha$ has not been injured by any $\gamma$
with $\alpha \subset \gamma \subseteq \beta$, $\phi_e(\langle
h(\alpha), d \rangle \downarrow)= w$, $[l, \infty ) \subseteq W_w$,
and $W_d \subseteq \bigsqcup_{\alpha \subset \delta \subseteq \beta}
F_\delta$.  Since $F_\delta = W_{g(\delta)}$, all this information is
$\Delta^0_3$ and hence can be easily coded into a tree. In the
interest of space and energy we are not going to go into the details of
the definition of the tree.  Similar constructions with all the
details can be found in Section~\ref{tree2}, \citet{Cholak:95},
\citet{Cholak:94*1}, and \citet{Weber:04}.  There is one added twisted
that there is finite injury along the true path.  But that kink was
discussed above and is implemented in the standard fashion.  \hfill
$\square$

\subsection{Reflecting on $\phi(A)$ and Theorem~\ref{orbit2}}
\label{restricts2}

Theorem~\ref{orbit2} implies that $\mathcal{O}$ is different than any
other known orbit.  But it might be worthwhile to reflect on
$\mathcal{O}$'s similarity to the orbit formed by the maximal sets or
the orbit formed by the Herrman sets (for a definition of Herrmann
sets, see \citet{Cholak.Downey.ea:01}).  This reflection will also
impact how we approach the proof of Theorem~\ref{orbit2}.

\begin{definition}\label{da}
  $\mathcal{D}(A)$ is the ideal generated by the sets $F$ such that
  either $F \cap A = \emptyset$ or $F \subseteq^* A$.
  $\mathcal{D}(A)$ is a $\Sigma^0_3$ ideal of $\E$.  Let
  $\E_{\mathcal{D}(A)}$ be $\E$ modulo $\mathcal{D}(A)$.  We write $X
  \subseteq_{\D(A)} Y$ if $X$ is contained in $Y$ modulo $\D(A)$. If
  $A$ is understood from the context we drop the ``$(A)$''.
\end{definition}

The last clause of $\phi(A)$ implies that $\E_{\mathcal{D}}$ is the
two element Boolean algebra.  This is also the case with maximal sets
and Herrmann sets. When this is the case we say that $A$ is
$\mathcal{D}$-maximal.  It is also possible to consider $A$ where
$\E_{\mathcal{D}}$ is a Boolean algebra, in which case, $A$ is
called $\mathcal{D}$-hhsimple.  (For more on $\mathcal{D}$-hhsimple
sets, see \citet{Cholak.Downey.ea:01}, \citet{Herrman.Kummer:94}, and
\citet{Kummer:91} in that order.)

Assume that $A$ is $\mathcal{D}$-hhsimple.  Furthermore assume that $W
\not=_{\mathcal{D}} A$. Then there is a $\tilde{W}$ such that
$W \cap \tilde{W} =_{\mathcal{D}} \emptyset$ and $W \cup
\tilde{W} =_{\mathcal{D}} \omega$.  So there is a set $F \in
\mathcal{D}$ such that $W \cap \tilde{W} \subseteq F$ and $W \cup
\tilde{W} \cup F = \omega$.  Therefore there is a computable set $R$
such that $R \cap \overline{F} = W \cap \overline{F}$.  

Let $\widetilde{\L}(A)$ be the definable (in $\E$) quotient
substructure of $\S_{\R}(A)$ given by $\{R \cap H: R \text{ is
  computable}\}$ modulo $\R(A)$.  Given the above paragraph, it is
straightforward to verify that $\widetilde{\L}(A)$ and
$\E_{\mathcal{D}}$ are $\Delta^0_3$-isomorphic.

Assume $A$ and $\Ahat$ are automorphic by $\Phi$.  By
Theorem~\ref{oldsplits}, $\S_{\R}(A)$ and $\S_{\R}(\Ahat)$ are
$\Delta^0_3$-isomorphic via an isomorphism induced by $\Phi$.  So
$\widetilde{\L}(A)$ and $\widetilde{\L}(\Ahat)$ are
$\Delta^0_3$-isomorphic via an isomorphism induced by $\Phi$.  Hence
$\E_{\mathcal{D}(A)}$ and $\E_{\mathcal{D}(\Ahat)}$ are
$\Delta^0_3$-isomorphic.  (A similar argument appeared in Section~11
of \citet{MR2004f:03077}.) Hence we have the following theorem.

\begin{theorem} \label{remove}
  Assume that $A$ is $\mathcal{D}$-hhsimple.  If $A$ and $\Ahat$ are
  automorphic via $\Phi$ then $\E_{\mathcal{D}(A)}$ and
  $\E_{\mathcal{D}(\Ahat)}$ are $\Delta^0_3$-isomorphic via an
  isomorphism induced by $\Phi$.
\end{theorem}

One should compare this theorem to Theorem~\ref{improved} where the
hypothesis that $A$ be $\mathcal{D}$-hhsimple is removed but the
complexity of the isomorphism increases to $\Delta^0_6$.

Soare showed that the maximal sets, $M$, do not form an effective
orbit by exploiting the fact that deciding if $W \subseteq^* M$ or
$W\cup M =^* \omega$ is $\Delta^0_3$. Soare built maximal sets $M$ and
$\widehat{M}$ such that for each computable function $f$ there is an
$e$ with $W_e \subseteq^* M$ iff $W_{f(e)} \cup \widehat{M} = ^*
\omega$.  (For more details, see \citet{Soare:74} and
\citet{Cholak:90}.)

But Theorem~\ref{remove} implies that we cannot exploit the fact of
deciding if $W \subseteq_{\mathcal{D}} A$ or $W =_{\mathcal{D}}
\omega$ is $\Delta^0_3$ to show there are $A$ and $\Ahat$ in
$\mathcal{O}$ which are not $\Delta^0_3$-isomorphic. Hence the
proposed approach of Cholak and Downey (thankfully unpublished) to the
proof of Theorem~\ref{orbit2} just cannot work.  To show
Theorem~\ref{orbit2} we exploited the fact that given a set $W$
disjoint from $A$ we cannot always computably in $\bf{0''}$ find an
$A$-special set disjoint from $W$.

\section{On the complexity of orbits of $\E$}
\label{restricts3}

The goal of this section is to improve Theorem~\ref{remove} and add to
our comments from Section~\ref{restricts2}.  We are going to do this
by coding where $W$, for $W \not=_{\mathcal{D}} A$, must go under an
arbitrary automorphism of $\E$, using various splits of $A$.  We will
break this into two subsections: the first subsection will focus on
the coding and the second subsection will present the results which use
this coding.

\subsection{Maximal supports} 

Fix a \ce set $A$.  A definition of $\mathcal{D}(A)$ can be found in
Definition~\ref{da}.

\begin{definition}
  $M$ is \emph{maximally supported} by $S$ if $M$ is supported by $S$
  (so $S$ is a split of $A$, $S \subseteq M$ and $(M-A) \sqcup S$ is a
  \ce set) and for all $W$, if $W$ is supported by $S$, then $W
  \subseteq_{\mathcal{D}} M \cup A$.
\end{definition}

\begin{lemma}
  Whether $S$ maximally supports $M$ is $\Pi^0_4$.
\end{lemma}

\begin{proof}
  By Lemma~\ref{supportsisSigma03}, whether $T$ supports $X$ is
  $\Sigma^0_3$.
\end{proof}

If $S$ is a maximal support of $W$ and $T =_{\R} S$ then $T$ is a
maximal support of $W$.

\begin{lemma} \label{msorder}
  If $Y \nsubseteq_{\mathcal{D}} X$, $S$ is a maximal support for
  $X$ and $T$ is a support for $Y$ then $T \nsubseteq_{\mathcal{\R}}
  S$.
\end{lemma}

\begin{proof}
  Since $S$ maximally supports $X$, $S$ cannot support $Y$.  So $T$ is
  not a subset of $S$.  The same holds modulo $\mathcal{R}(A)$.
\end{proof}

Note it is possible that $S$ and $T$ maximally support $W$ but $S
\neq_{\R} T$. But this will not cause a problem.

Recall $A$ is promptly simple iff there is a computable function $p$
such that for all $W$, if $W$ is infinite, then there is an $x$ and $s$
such that $x \in W_{\text{at }s} \cap A_{p(s)}$.  Also if $A$ is
simple then $W \subseteq_{\mathcal{D}} M$ iff $W \subseteq^* M \cup
A$.

\begin{lemma}\label{pmsupport}
  Assume that $A$ is promptly simple.  Let $A \subseteq M$.  There is
  an $S$ such that $M$ is maximally supported by $S$.
  
  Furthermore $S = M \searrow A$ using $\{A_{p(s)}\}_{s \in \omega}$
  as the enumeration of $A$; i.e., $S$ is the set of $x$ such that $x$
  enters $M$ at stage $s$ and $x$ is not in $A_{p(s)}$ but $x$ is in
  $A$.
\end{lemma}

\begin{proof}
  $M$ is supported by the $S$ defined above; $(M-A) \sqcup S$ is the
  set of $x$ such that $x$ enters $M$ at stage $s$ and $x$ is not in
  $A_{p(s)}$.
  
  To ensure $M$ is maximally supported by $S$ it is enough to show the
  following conditions are met:
  \begin{equation*}
    \tag*{$\mathcal{N}_{e,i}$:} \text{either } W_e \subseteq^* 
    M \cup A \text{ or } W_i \neq (W_e - A) \sqcup S.
  \end{equation*} 
  Assume $W_e \not\subseteq^* M \cup A$ and $W_i = (W_e - A) \sqcup S$
  (i.e., that we fail to meet $\mathcal{N}_{e,i}$). Then $W = (W_e
  \cap W_i) \backslash (M \cup A)$ is infinite.  Then there is an $x$
  and $s$ such that $x \in W_{\text{at }s} \cap A_{p(s)}$. Now $x$ is
  in $W_i$ and thus in one of $W_e - A$ or $S$.  But $x$ cannot be in
  either of these two sets.  Contradiction.
\end{proof}

It would be nice if we could prove the above lemma for all $A$ but the
above proof heavily relies on the assumption that $A$ was promptly
simple.  However we do have the following lemma.

\begin{lemma}\label{msupport}
  For all $W$, $\tilde{W}$, if $W \neq_{\mathcal{D}} \tilde{W}$ then
  there is an $M$ such that $M$ is maximally supported by $S= M
  \searrow A$, $M \subseteq W$, and $M \nsubseteq_{\mathcal{D}}
  \tilde{W}$.
\end{lemma}

\begin{proof}
  Fix $W$ and $\tilde{W}$. Clearly $M \searrow S$ supports $M$.  So we
  must build $M$ to meet the following requirements:
 \begin{equation*}
    \tag*{$\mathcal{N}_{e,i}$:} \text{either } W_e 
   \subseteq_{\mathcal{D}(A)} 
    M \cup A \text{ or } W_i \neq (W_e - A) \sqcup S.
  \end{equation*} 
  (i.e., either $W_e$ is contained in $M \cup A$ modulo $\D(A)$ or $S$
  does not support it) and
  \begin{equation*}
    \tag*{$\mathcal{P}_{e,i}$:} \text{either } W_e 
   \cap A \neq \emptyset, \text{ or }   W_i 
   \cap A \neq \emptyset, \text{ or } M \cup W_i \cup A \nsubseteq 
   \tilde{W} \cup W_e \cup A
  \end{equation*} 
  (so $M$ is not contained modulo $\D(A)$ in $\tilde{W}$).  Assume
  that these requirements are linearly ordered.
  
  To meet $\mathcal{N}_{e,i}$ we will hold everything in $X = (W_e \cap
  W_i) \backslash (M \cup A)$ out of $M$ until there is an $x \in X
  \cap A$ and hence $W_i \neq (W_e - A) \sqcup S$.  Assume this fails.
  Then $X$ is disjoint from $A$.  So if $W_i = (W_e - A) \sqcup S$
  then $W_e \subseteq M \cup A \cup X$. And hence we still meet
  $\mathcal{N}_{e,i}$.
 
  To meet $\mathcal{P}_{e,i}$ we need to first define a length of
  agreement function (to measure a $\Pi^0_2$ fact).  Let $l(s)$ be the
  greatest $x$ such that $(W_{e,s} \cap A_s) \restriction x =
  \emptyset$, $(M_s \cup W_{i,s} \cup A_s) \restriction x =
  (\tilde{W}_s \cup W_{e,s} \cup A_s) \restriction x$, and $(W_{i,s}
  \cap A_s) \restriction x =\emptyset$.  Let $m(0)=0$. If $l(s) >
  m(s-1)$ then $s$ is \emph{expansionary} (for $\mathcal{P}_{e,i}$)
  and $m(s)=l(s)$; otherwise $m(s)=m(s-1)$.
  
  If there are infinitely many expansionary stages we must take some
  action to ensure $\mathcal{P}_{e,i}$ is met. At expansionary stages
  we will dump everything in $W$ which is not restricted by higher
  priority requirements into $M$ and reset all lower priority
  requirements.
  
  As we argued above, the set $X$ of $x$ which is restrained by
  high priority requirements is disjoint from $A$.  Therefore if there
  are infinitely many expansionary stages then $M \cup Z \cup A= W
  \cup Z \cup A$, where $Z$ is the union of finitely many $X$s from the
  higher priority negative requirements.  Hence $W =_{\mathcal{D}} M
  =_{\mathcal{D}} \tilde{W}$.  Hence, under the above hypothesis,
  there cannot be infinitely many expansionary stages and
  $\mathcal{P}_{e,i}$ is met.
\end{proof}

\subsection{Coding with maximal supports}

\begin{theorem}\label{orbitsps}
  Assume that $A$ and $\Ahat$ are promptly simple.  Then $A$ and
  $\Ahat$ are automorphic iff $A$ and $\Ahat$ are $\Delta^0_3$
  automorphic.
\end{theorem}

\begin{proof}
  Assume that $A$ and $\Ahat$ are automorphic via $\Phi$. We can
  assume that $\Phi \restriction \E^*(A)$ is $\Delta^0_3$. We must
  show that $\Phi \restriction \mathcal{L}^*(A)$ is $\Delta^0_3$.  We
  know that $\mathcal{S}_{\R}(A)$ and $\mathcal{S}_{\R}(\Ahat)$ are
  $\Delta^0_3$ isomorphic via an isomorphism $\Theta$ induced by
  $\Phi$.
  
  Given $W$, look for a support $S$ of $W$, a set $\tilde{W} \subseteq
  \widehat{\omega}$, and a support $\tilde{S}$ of $\tilde{W}$ such
  that $ S \subseteq_{\mathcal{R}} \Theta^{-1}(\tilde{W} \searrow
  \Ahat)$ and $ \tilde{S} \subseteq_{\mathcal{R}} \Theta(W \searrow
  A)$.  Such sets exist; consider $\tilde{W} = \Phi(W)$, $S =
  \Phi^{-1}(\Phi(W) \searrow \Ahat)$, and $\tilde{S} = \Phi(W \searrow
  A)$.  Since such sets exist, we can find them using $\bf{0''}$ as
  an oracle.
  
  Since $\Theta$ is induced by the automorphism $\Phi$, by
  Lemma~\ref{pmsupport}, $\Theta^{-1}(\tilde{W} \searrow \Ahat)$
  maximally supports $\Phi^{-1}(\tilde{W})$.  Therefore, by
  Lemma~\ref{msorder} and the fact that for simple sets, $A$, $=^*$,
  and $=_{\mathcal{D}}$ agree, $W \subseteq^* \Phi^{-1}(\tilde{W})$.
  Similarly $\tilde{W} \subseteq^*\Phi(W) $.  So $W =^*
  \Phi^{-1}(\tilde{W})$ and $\tilde{W} =^* \Phi(W)$ and hence
  $\tilde{W} =^* \Phi(W)$.
\end{proof}

\begin{theorem}\label{improved}
  If $A$ and $\Ahat$ are automorphic via $\Phi$ then
  $\E_{\mathcal{D}(A)}$ and $\E_{\mathcal{D}(\Ahat)}$ are
  $\Delta^0_6$-isomorphic via an isomorphism induced by $\Phi$.
\end{theorem}

\begin{proof}
  Assume that $A$ and $\Ahat$ are automorphic via $\Phi$. We can
  assume that $\Phi \restriction \E^*(A)$ is $\Delta^0_3$.  We know
  that $\mathcal{S}_{\R}(A)$ and $\mathcal{S}_{\R}(\Ahat)$ are
  $\Delta^0_3$ isomorphic via an isomorphism $\Theta$ induced by
  $\Phi$.  Given $W$ we must find a $\tilde{W}$, in a $\Delta^0_6$
  way, such that $\tilde{W} =_{\mathcal{D}} \Phi(W)$.
  
  By Lemma~\ref{msupport}, $Y \subseteq_{\mathcal{D}} \tilde{Y}$ iff,
  for all $M$ and $X$, if $M \subseteq Y$, $M$ is maximally supported
  by $S = M \searrow A$, and $S$ supports $X$, then $X
  \subseteq_{\mathcal{D}} \tilde{Y}$.  Since $\Theta$ is induced by
  the automorphism $\Phi$, $\tilde{W} \subseteq_{\D} \Phi(W)$ iff for
  all $\tilde{M}$ and $X$, if $\tilde{M} \subseteq \tilde{W}$,
  $\tilde{M}$ is maximally supported by $\tilde{S} = \tilde{M}
  \searrow \Ahat$, and $\Theta^{-1}(\tilde{S})$ supports $X$, then $X
  \subseteq_{\mathcal{D}} W$, a $\Pi^0_5$-statement.  And similarly,
  $\Phi(W) \subseteq_{\D} \tilde{W}$ iff for all $M$ and $\tilde{X}$,
  if $M \subseteq W$, $M$ is maximally supported by $S = M \searrow A$,
  and $\Theta(S)$ supports $\tilde{X}$, then $\tilde{X}
  \subseteq_{\mathcal{D}} \tilde{W}$, a $\Pi^0_5$-statement.
  
  Therefore whether $\tilde{W} =_{\D} \Phi(W)$ is $\Pi^0_5$.  Since
  such a $\tilde{W}$ exists, it can be found using $\bf{0^{(5)}}$ as
  an oracle.
\end{proof}

\begin{corollary}
  If $A$ is simple, then $A$ and $\Ahat$ are automorphic iff $A$ and
  $\Ahat$ are $\Delta^0_6$-automorphic.
\end{corollary}

\begin{proof}
  Assume that $A$ and $\Ahat$ are automorphic by $\Phi$ where $\Phi
  \restriction \E^*(A)$ is $\Delta^0_3$.  Since $A$ is simple, if $W
  \subseteq \overline{A}$ then $W$ is finite.  Therefore
  $\mathcal{L}^*(A)$ and $\E_{\mathcal{D}(A)}$ are isomorphic, by the
  identity map.  Therefore $\Phi \restriction \L^*(A)$ is
  $\Delta^0_6$. So $\Phi$ is $\Delta^0_6$.
\end{proof}

If $A$ is simple and $A \subset W$ then where an automorphism of $\E$
takes $W$ is completely determined by certain splits of $A$, the
maximal supports. Hence the following is a corollary of the proofs of
Theorem~\ref{improved} (\ref{orbitsps}) and Theorem~\ref{orbits}.

\begin{theorem}\label{simpleo}
  The (promptly) simple sets $A$ and $\Ahat$ are automorphic iff there
  are $\Psi$, $\B$, $\widehat{\B}$, and $\Theta$ such that
  \begin{enumerate}
  \item $\L^*(A)$ and $\L^*(\Ahat)$ are $\Delta^0_6$-isomorphic
    ($\Delta^0_3$-isomorphic) via $\Psi$,
  \item $\B$ and $\widehat{\B}$ are extendible algebras which are
    extendibly $\Delta^0_3$ isomorphic via $\Theta$,
  \item $\B$ supports $\L^*(A)$,
  \item $\widehat{\B}$  supports $\L^*(\Ahat)$,
  \item the isomorphisms $\Psi$ and $\Theta$ preserve supports.
 \end{enumerate}
\end{theorem}

The $r$-maximal sets are simple. So $r$-maximal sets are automorphic
iff they are $\Delta^0_6$-automorphic.  But this is not a ``nice''
algebraic classification, at least for $r$-maximal sets.  It is
possible that the $\L^*$s of $r$-maximal sets have a nice structure.
So we might be able to replace Condition~1 of Theorem~\ref{simpleo}
with something more algebraic and easier to understand, like the other
conditions.  The reader is directed to the last section of
\citet{Cholak.Nies:99} for some suggestions.  We should point out that
\citet{Lempp.Nies.ea:01} have shown that there is no $\Delta^0_3$
classification (``nice'' or otherwise) of the $\L^*$s of $r$-maximal
sets.  But this does not rule out a ``nice'' arithmetic classification
of the $\L^*$s.

The results in this section and that of Section~\ref{restricts2}
drive home the point that to build sets whose orbits are complex we
are forced to use techniques like those described in
Sections~\ref{sec:orbit-mathcalo-4} and \ref{sec:changes-needed-proof}.
In a forthcoming paper we will do just that.

\end{document}